\input amstex
\documentstyle{amsppt}
\magnification=\magstep1
\hsize= 6 true in
\hoffset= 0.25 true in
\vsize= 9 true in
\NoRunningHeads
\define\R{\bold R}
\define\RR{\bold R^2}
\define\RRR{\bold R^3}
\define\I{[0,1]}

\define\p{^{\prime}}
\define\pp{^{\prime\prime}}
\define\bproper{$\partial$-proper}
\define\bparallel{$\partial$-parallel}
\define\bd{\partial}
\define\halfspace{\bold R^2 \times [0,1)}
\define\prodspace{\bold R^2 \times [0,1]}
\define\homeo{homeomorphic}
\define\3m{3-manifold}
\define\inc{incompressible}
\define\irr{irreducible}
\define\birr{$\partial$-irreducible}
\define\binc{$\partial$-incompressible}
\define\ann{anannular}
\define\nzsp{$n_0$-standard position}
\define\nzmp{$n_0$-monotone position}
\define\inte{int \, }
\define\Inte{Int \, }
\define\FF0{{\Cal F}_0}
\define\Fn0{F_{n_0}}
\define\FFn0{\Cal F_{n_0}}
\define\Cn0{C_{n_0}}
\define\Yn0{Y_{n_0}}
\define\PP{\Cal P}
\define\QQ{\Cal Q}
\define\JJ{\Cal J}
\define\AAA{\Cal A}
\define\al{\alpha}
\define\be{\beta}
\define\ga{\gamma}
\define\de{\delta}
\define\la{\lambda}
\define\La{\Lambda}
\define\si{\sigma}
\define\Si{\Sigma}
\define\th{\theta}
\define\thh{\widehat{\theta}}
\define\Bh{\widehat{B}}
\define\Ch{\widehat{C}}
\define\Yt{\widehat{Y}}

\topmatter
\title Attaching Boundary Planes to Irreducible Open 3-Manifolds \endtitle 
\author Robert Myers \endauthor
\address Department of Mathematics, Oklahoma State University, 
Stillwater, OK 74078-0613, USA  \endaddress
\email myersr\@math.okstate.edu \endemail

\thanks This research was partially supported by a summer research award 
from the College of Arts and Sciences at Oklahoma State University. 
\endthanks

\abstract 
Given any connected, open 3-manifold $U$ having finitely 
many ends, a non-compact 3-manifold $M$ is constructed having the 
following properties: 
the interior of $M$ is homeomorphic to $U$; 
the boundary of $M$ is the disjoint union of finitely many planes; 
$M$ is not almost compact; $M$ is eventually end-irreducible;  
there are no proper, incompressible embeddings of $S^1 \times \bold R$ 
in $M$; 
every compact subset of $M$ is contained in a larger compact subset whose 
complement is anannular; 
there is a compact subset of $M$ whose complement is 
$\bold P^2$-irreducible. 

If $U$ is irreducible it also has the following two properties: 
every proper, non-trivial plane in $M$ is boundary-parallel; 
every proper surface in $M$ 
each component of which has non-empty boundary and is non-compact 
and simply connected lies in a collar on $\partial M$. 

This construction 
can be chosen so that $M$ admits no homeomorphisms which take 
one boundary plane to another or reverse orientation. For the given $U$ 
there are uncountably many non-homeomorphic such $M$. 

Two auxiliary results may be of independent interest. First, general 
conditions are given under which infinitely many ``trivial'' compact 
components of the intersection of two proper, non-compact surfaces in 
an irreducible 3-manifold can be removed by an ambient isotopy. 
Second, $n$ component tangles in a 3-ball are constructed such that 
every non-empty union of components of the tangle has hyperbolic 
exterior. 

\endabstract 

\endtopmatter
\document

\pagebreak

\head Introduction  \endhead

Suppose $U$ is an open 3-manifold and $M$ is a non-compact 3-manifold 
each of whose boundary components is \homeo\  to $\RR$. 
If $U$ is \homeo\  to the interior of $M$, then one may say that 
$M$ is obtained by {\bf attaching boundary planes} to $U$. The fact that 
this can be done in different ways is most dramatically illustrated by 
the existence of 3-manifolds with interior \homeo\  to $\RRR$ 
and boundary \homeo\  to $\RR$ which are not \homeo\  to 
$\RR \times [0,\infty)$. The first such example was constructed by 
Fox and Artin \cite{9}. Tucker \cite{18} later gave a different method for constructing 
such examples and showed that planes were the only surfaces which 
could be, in a sense explained below, ``bad'' boundary components of 
3-manifolds.  

A non-compact 3-manifold $M$ is {\bf almost compact} if there is a 
compact 3-manifold $Q$ and a closed subset $K$ of $\bd Q$ such that 
$M$ is \homeo\ to $Q-K$. The examples mentioned above fail to be almost 
compact. Tucker \cite{18} showed that if one attaches a connected surface 
with finitely generated fundamental group to a $\bold P^2$-\irr, almost compact 
3-manifold in such a way that it becomes an \inc\ boundary component, then 
either the new 3-manifold is almost compact or the surface is a plane. 
Scott and Tucker \cite{17} gave other examples of 3-manifolds 
which are not almost compact but have almost compact interiors, 
including an example with boundary consisting of two disjoint planes 
such that the complement of each plane is \homeo\  to $\halfspace$ and 
an example with boundary a single plane whose complement is \homeo\  to 
$S^1 \times \RR$. It can be shown that the latter example contains no 
proper non-separating planes. 

This paper gives a general procedure for attaching boundary planes to $U$ 
to obtain an $M$ having special properties which may be very 
different from those of $U$. The open 3-manifold 
$U$ is required to be connected and irreducible and to have 
finitely many ends. It is not assumed to be orientable or $\bold P^2$-\irr. 
One attaches any finite number of boundary planes to $U$ 
subject to the restriction that one attaches at least one plane to each 
end. The resulting 3-manifold $M$ is not almost compact. However, it has 
an important property which is shared by those almost compact 3-manifolds 
whose boundaries are finite disjoint unions of planes, namely eventual 
end-irreducibility.  This property, introduced by Brown \cite{8}, ensures that the 
ends of the 3-manifold, while perhaps not tame, are at least not excessively wild, 
in the sense that they can be analyzed using \inc\ surface theory. 
This has proven to be a fruitful concept in the study of non-compact 3-manifolds. 
See \cite{2, 3, 4, 6, 7, 8}.  
The next section includes discussions of ends and of 
eventual end-irreducibility. 

The most important special properties considered in this paper concern 
certain embeddings of surfaces. A surface $S$ embedded in a 3-manifold 
$M$ is {\bf proper} if $S \cap \bd M=\bd S$ and $S \cap C$ is compact 
for every compact subset $C$ of $M$. A proper surface $S$ with 
$\bd S = \emptyset$ is {\bf \bparallel\ }  if some component of $M-S$ has 
closure \homeo\  to $S \times [0,1]$ with $S=S \times \{0\}$ and 
$S \times \{1\}$ a component of $\bd M$; it is {\bf end-parallel} or 
{\bf trivial} if some component of $M-S$ has closure \homeo\  to 
$S \times [0,\infty)$ with $S=S \times \{0\}$. A 3-manifold $M$ is 
{\bf aplanar} if every proper plane in $M$ is either \bparallel\  
or end-parallel; it is {\bf acylindrical} if the same is true for every 
proper incompressible cylinder $S^1 \times {\bold R}$. It is {\bf totally 
acylindrical} if it contains no proper \inc\ cylinders. 
The boundary planes will be attached so that $M$ is 
aplanar and totally acylindrical. 
In particular one can take an irreducible one ended open 3-manifold which 
contains non-trivial planes or non-trivial \inc\ cylinders, such as the interior of a cube with 
handles or $S^3-K$ for $K$ a torus knot, cable knot, or composite knot, 
and create an aplanar, totally acylindrical 3-manifold by attaching a 
single boundary plane. 

It will be shown that $M$ has two further embedding properties. It is 
{\bf strongly aplanar}, meaning that in addition to being aplanar  
it has the property that given any proper surface $\Cal P$ such that each 
component of $\Cal P$ is non-compact, is simply connected, and has non-empty 
boundary, there exists a collar on $\bd M$ containing $\Cal P$. It is also 
{\bf anannular at infinity} in the sense that for every compact subset 
$K$ of $M$ there is a compact subset $L$ of $M$ containing $K$ such that 
$M-L$ is anannular. Note for comparison that $\halfspace$ is aplanar but not 
strongly aplanar, while $\prodspace$ is aplanar but neither strongly 
aplanar nor anannular at infinity. (See \cite{15}.) These two properties are involved 
in the study of ``plane sums'' of non-compact 3-manifolds, i.e\. 
3-manifolds obtained by gluing together a collection of 3-manifolds 
whose boundary components are planes along these planes. In \cite{15} it will 
be shown that (subject to some mild additional hypotheses) if the 
summands are irreducible, strongly aplanar, and anannular at infinity, 
then the image of each gluing plane is non-trivial in the sum, and every 
non-trivial plane in the sum is ambient isotopic to one of these planes. 
This result is then used to investigate a non-compact analogue of the 
connected sum called the ``end sum.'' The present paper provides the 
mechanism for generating the relevant examples. 

By a modification of the basic construction $M$ can be built so that,  
in addition to the previous properties, it admits no orientation reversing 
self-homeomorphisms, there are no self-homeomorphisms taking one boundary 
plane to another, and there are uncountably many pairwise 
non-homeomorphic such $M$ having the same number of boundary planes 
per end of $U$. Moreover all these properties hold for the 3-manifolds obtained 
by deleting any collection of boundary planes from $M$ subject to the 
restriction that there remains at least one boundary plane attached to each 
end. These properties are also relevant to the study of 
plane sums and end sums. 

Although \irr\ 3-manifolds are the main objects of interest in this paper 
it is worth noting that some of our results generalize to the reducible 
case. In fact the only properties of $M$ which require the irreducibility 
of $U$ are aplanarity and strong aplanarity. Moreover $M$ can be 
constructed so as to be {\bf eventually $\bold P^2$-\irr} in the sense that 
there is a compact subset of $M$ whose complement is $\bold P^2$-\irr. 
One can, for example, create a totally acylindrical, 
eventually $\bold P^2$-\irr\ 3-manifold by attaching two boundary planes 
to the product of a closed, connected surface with $\R$, including the 
cases when the surface is $S^2$ or $\bold P^2$.

Two auxiliary results in this paper may be of independent interest. 
First, we give general conditions under which infinite collections of  
``trivial'' intersection curves of two non-compact proper surfaces in an irreducible 
3-manifold can be removed by an ambient isotopy. By a ``trivial'' intersection 
curve we mean a simple closed curve which bounds disks on both surfaces 
or a proper arc which is \bparallel\ on both surfaces. These results are used 
both in the present paper and in \cite{15}. 
Second, we prove  
the existence of ``poly-excellent tangles.''  A poly-excellent $n$-tangle is 
the union of $n$ disjoint proper arcs in a 3-ball such that the union of any 
non-empty collection of its components has hyperbolic exterior. 
This result is required in our modification of the basic construction. 
 
The paper is organized as follows. Section 1 contains background material 
and discusses exhaustions of non-compact 3-manifolds. In particular it 
introduces the concept of a ``nice'' exhaustion. It is readily seen that a
3-manifold $M$ with a nice exhaustion is not almost compact but is 
eventually end-irreducible, eventually ${\bold P}^2$-\irr, and \ann\ at 
infinity. Section 2 gives conditions under which 
``trivial'' intersections of non-compact surfaces can be removed by 
an ambient isotopy. Section 3 reformulates some work of Winters \cite{19, 20} 
to show that a \3m\  with a nice exhaustion is totally acylindrical and, if it 
is \irr, is aplanar. Sections 4 and 5 show that an \irr\ \3m\ with a nice 
exhaustion is strongly aplanar.  Section 6  
shows how to construct $M$ from $U$ so that $M$ has a nice exhaustion; 
it also describes the modification of  this basic construction and proves 
the additional properties listed above. The proof of the existence of 
poly-excellent tangles has a different flavor from the rest of the paper and 
is given in an appendix so as not to disrupt the main line of the argument.

\head 1. Preliminaries  \endhead

We shall work throughout in the PL category. An $m$-manifold $M$ may or 
may not have boundary but is assumed to be second countable. 
$\bd M$ and $\inte M$ denote the manifold theoretic boundary and interior  
of $M$, respectively. Let $A$ be a subset of $M$. The topological boundary, 
interior, and closure of $A$ in $M$ are denoted by $Fr_M A$, $Int_M A$, and 
$Cl_M A$, respectively, with the subscript deleted when $M$ is clear 
from the context. All isotopies of $A$ in $M$ will be ambient. 
$A$ is {\bf bounded} if $Cl \, A$ is compact. $M$ is 
{\bf open} if $\bd M=\emptyset$ and no component of $M$ is compact. 

A {\bf  surface} is a 2-manifold; no assumptions are made about its being 
connected or compact or having a boundary.

A map $f:M \rightarrow N$ of manifolds is {\bf \bproper\ } if 
$f^{-1}(\bd N)=\bd M$. It is {\bf end-proper} if preimages of compact 
sets are compact. It is {\bf proper} if it has both these properties. 
These terms are applied to a submanifold if its inclusion map has the 
corresponding property. 

Let $S$ be a proper codimension one submanifold of the $m$-manifold $M$. 
Suppose $S\p$ is either another such submanifold such that $\inte S \cap 
\inte S\p=\emptyset$ or is an end-proper submanifold of 
$\bd M$. Assume that $\bd S = \bd S\p$. Then $S$ and $S\p$ are 
{\bf parallel} if some component of $M-(S \cup S\p)$ has closure \homeo\ 
to $S \times \I$ with $S \times \{0\}=S$ and $S \times \{1\}=S\p$ 
when $\bd S = \emptyset$, while $((\bd S) \times \I) \cup (S \times \{1\}) 
=S\p$ when $\bd S \neq \emptyset$. The product $S \times \I$ is a 
{\bf parallelism} between $S$ and $S\p$. When $S\p \subseteq \bd M$ 
one says that $S$ is {\bf \bparallel\ }. We say that $S$ is 
{\bf end-parallel} or {\bf trivial} if some component of $M-S$ has 
closure \homeo\ to $S \times [0,\infty)$ with $S=S \times \{0\}$. 

Infinite sequences, unless indicated otherwise, will be indexed by 
the set of non-negative integers. 

An {\bf exhausting sequence} $C=\{C_n\}$ for a non-compact $m$-manifold 
$M$ is a sequence 
$C_0 \subseteq C_1 \subseteq C_2 \subseteq \cdots$ of compact subsets 
of $M$ whose union is $M$. A sequence 
$V_0 \supseteq V_1 \supseteq V_2 \supseteq \cdots$ of open subsets of 
$M$ is an {\bf end sequence} associated to $C$ if each $V_n$ is 
a component of $M-C_n$. Two end sequences $\{V_n\}$ and $\{W_p\}$ 
associated to exhausting sequences $C$ and $K$ for $M$ 
are {\bf cofinal} if for every $n$ there is a $p$ such that 
$V_n \supseteq W_p$ and for every $p$ there is a $q$ such that 
$W_p \supseteq V_q$. Cofinality is an equivalence relation on 
end sequences of $M$. The equivalence classes are called the 
{\bf ends} of $M$. The set of all ends of $M$ is denoted by 
$\varepsilon (M)$. An end-proper map $M \rightarrow N$ induces a 
well defined function $\varepsilon(M) \rightarrow \varepsilon(N)$. 
If $\bd M$ has no compact components, then the inclusion map induces 
a well defined bijection $\varepsilon(int \, M) \rightarrow \varepsilon(M)$. 

An exhausting sequence $C$ for a connected, non-compact  $m$-manifold $M$ is an 
{\bf exhaustion} for $M$ if each $C_n$ is a compact, connected $m$-manifold,  
$C_n\cap\bd M$ is either empty or an $(m-1)$-manifold, 
$C_n \subseteq Int \, C_{n+1}$, and $M-C_n$ has no bounded 
components. Connected non-compact $m$-manifolds always have exhaustions. 
Given an exhaustion $C$ for $M$ and a subsequence $\{n_k\}$ of the 
non-negative integers, let $C\p_k=C_{n_k}$. Then $C\p$ is also 
an exhaustion for $M$ and will be called a {\bf subexhaustion} of $C$. 

The reader is referred to \cite{10} or \cite{11} for basic \3m\ topology, including 
the definition of incompressible surface. We adopt 
the conventions of \cite{11} that every proper disk in a \3m\ $M$ is \inc\  
and that a proper 2-sphere is compressible if and only if it bounds a 
$3$-ball. $M$ is {\bf \irr\ } if every 2-sphere in $M$ is compressible; 
it is $\bold P^2${\bf -\irr} if it contains no 2-sided projective planes. 
It is {\bf \birr\ } if $\bd M$ is \inc\ in $M$. It is {\bf \ann\ } 
if every proper \inc\ annulus in $M$ is \bparallel\. It is  
{\bf atoroidal} if every proper \inc\ torus in $M$ is \bparallel\. 

A {\bf partial disk} is a pair $(D,\bd_0 D)$, where 
$D$ is a disk and $\bd_0 D$ is a non-empty finite union of disjoint arcs in 
$\bd D$; the {\bf order} of $(D,\bd_0 D)$ is the number of these arcs. 
A {\bf halfdisk} is a partial disk of order one;  
a {\bf band} is a partial disk of order two. 
$\bd D -\inte \bd_0 D$ is denoted by $\bd_1 D$. A partial disk may be 
denoted by $D$ when $\bd_0 D$ is clear from the context. 
Suppose $D$ is a partial disk contained in a surface $S$. Then $D$ 
is {\bf proper} in $S$ if $D \cap \bd S = \bd_0 D$. If $D$ is a proper 
partial disk in $S$ such that no component of $S-Int_S D$ is a proper 
halfdisk $D\p$ in $S$ with $\bd_1 D\p=D \cap D\p$, then $D$ is 
{\bf well embedded} in $S$.

A proper surface $S$ in $M$ is {\bf \binc\ } if it is not a \bparallel\ 
disk and whenever $D$ is a halfdisk in $M$ such that $D \cap \bd M     
=\bd_0 D$ and $D \cap S = \bd_1 D$, one has that $\bd_1 D$ is 
\bparallel\  in $S$. 

We now give some terminology for some standard isotopies which will 
be used later. Suppose $S$ and $T$ are end-proper surfaces in an 
\irr\ \3m\ $M$. Suppose $S$ and $T$ are in general position and $J$ 
is a simple closed curve component of $\inte S \cap \inte T$ which bounds a disk 
$D$ on $S$ and a disk $G$ on $T$. Then $J$ is {\bf innermost} on $S$ if 
$D \cap T=J$. In this case there is a 3-ball $B$ in $M$ with 
$\bd B=D \cup G$. Let $B^+$ be a regular neighborhood of $B$ in $M$. 
There is an ambient isotopy of $S$ in $M$ supported in $B^+$ which 
carries $D$ to $G$ and then off $G$ into $B^+-B$. This is called a 
{\bf disk push} of $D$ across $B$ past $G$. Let $S^*$ and $D^*$ be the 
images of $S$ and $D$, respectively, under this isotopy. Then $S^*$ 
is in general position with respect to $T$ and 
$((S^*-D^*) \cap T) \subseteq ((S-D) \cap T)$ and $D^* \cap T=\emptyset$. 

Now suppose that $S$ and $T$ are also \bproper\  in $M$ and that 
there is an end-proper surface $R$ in $\bd M$ such that 
$\bd S \cup \bd T$ lies in $\inte R$ and $R$ is \inc\ in $M$.  
Suppose $\al$ is a component of $S \cap T$ which is an arc such that $\al=
\bd_1 D = \bd_1 G$, where $D$ and $G$ are proper halfdisks in $S$ 
and $T$, respectively. Then $\al$ is {\bf innermost} on $S$ if 
$D \cap T=\al$. In this case $\bd_0 D \cup \bd_0 G = \bd D\p$ 
for a disk $D\p$ in $\inte R$, and there is a 3-ball $B$ in $M$ 
such that $\bd B=D \cup G \cup D\p$. Let $B^+$ be a regular 
neighborhood of $B$ in $M$. There is an ambient isotopy of $S$ 
in $M$ supported in $B^+$ which carries $D$ to $G$ and then off 
$G$ into $B^+-B$. This is called a {\bf halfdisk push} of $D$ across $B$ 
past $G$. The images $S^*$ and $D^*$ then satisfy the same conditions 
as for a disk push. 

A {\bf partial plane} $P$ is a non-compact simply connected 2-manifold with 
$\bd P \neq \emptyset$. When $\bd P$ has exactly one component $P$ is called 
a {\bf halfplane}. We next give criteria for a proper plane or halfplane 
to be trivial. Note that a proper halfplane is trivial if and only if 
it is \bparallel\. 

\proclaim{Lemma 1.1} Let $M$ be a connected, \irr, non-compact \3m\. 
\roster
\item A proper plane $P$ in $M$ is trivial if and only if there 
exist sequences $\{D_n\}$ and $\{D_n\p\}$ of disks in $M$ such that $\{D_n\}$  
is an exhaustion for $P$, $D_n\p \cap P=\bd D_n$, and $\cup D_n\p$ is 
end-proper in $M$. 
\item A proper halfplane $P$ in $M$ is trivial if and only if there 
exist sequences $\{D_n\}$ and $\{D_n\p\}$ of halfdisks in $M$ such that 
$D_n$ is proper in $P$, $\{D_n\}$ is an exhaustion for $P$, $D_n\p \cap 
P=\bd_1 D_n=\bd_1 D_n\p$, $D_n\p \cap \bd M=\bd_0 D_n\p$, 
$\cup D_n\p$ is end-proper in $M$, and $\bd_0 D_n \cup \bd_0 D_n\p$ bounds 
a disk in $\bd M$. 
\endroster \endproclaim

\demo{Proof} Necessity is obvious in both cases. 

(1) Given any compact subset $K$ of $M$ there is an $n$ such that $D_n$ 
contains $K \cap P$ and $D_m\p \cap K=\emptyset$ for all $m \geq n$. 
If $K$ is a simple closed curve which meets $P$ transversely in a 
single point, then $K$ meets the 2-sphere $D_n \cup D_n\p$ transversely 
in a single point, contradicting the fact that $D_n \cup D_n\p$ bounds 
a 3-ball in $M$. Thus $P$ must separate $M$ into two components with 
closures $X$ and $Y$. By passing to a subsequence we may assume that all 
$D_n\p$ lie in $X$. Then for each compact subset $K$ of $X$ there is an $n$ 
such that $K$ lies in the 3-ball in $X$ bounded by $D_n \cup D_n\p$. 
It follows that $X$ is \homeo\  to $\halfspace$. 

(2) The proof is similar to that of (1) and is left to the reader. 
\qed \enddemo 
 
Let $M$ be a connected, non-compact \3m\. 
It is said to be {\bf eventually end-irreducible} if it 
contains a compact subset $J$ such that for every compact subset $K$ 
containing $J$ there is a compact subset $L$ containing $K$ such that 
every loop in $M-L$ which is null-homotopic in $M-J$ must be null-homotopic 
in $M-K$. 

Let $C$ be an exhaustion for a connected, non-compact \3m\ $M$. 
Denote $Fr \, C_n$ by $F_n$, $C_{n+1}-Int \, C_n$ by $X_{n+1}$, and 
$\cup_{n \geq k}F_n$ by $\Cal F_k$. The exhaustion is {\bf good} if for 
each $n \geq 0$ one has that $F_n \cup F_{n+1}$ is incompressible in 
$X_{n+1}$, 
$X_{n+1}$ is $\bold P^2$-irreducible, no component of $\FF0$ is a disk, and no  
component of $X_{n+1}$ has the form $F \times \I$, where $F \times \{0\}$    
and $F \times \{1\}$ are components of $F_n$ and $F_{n+1}$ respectively.       
Standard arguments show that in this case $\Cal F_0$ is \inc\ in $M-Int \, C_0$ 
and has no 2-sphere  
or projective plane components and that $M-Int \, C_0$ is $\bold P^2$-\irr\. 
It is easily seen that every subexhaustion of a good exhaustion is good. 
Note that $M$ itself need not be $\bold P^2$-\irr\ or even \irr. 

\proclaim{Lemma 1.2} If the connected, non-compact \3m\ $M$ 
admits a good exhaustion, then $M$ is eventually end-\irr\ and eventually 
$\bold P^2$-\irr\ and is not almost compact. 
\endproclaim

\demo{Proof} The first property is well known; its proof will be 
sketched for completeness. 
Let $J=C_0$. Suppose $K$ is compact and contains $J$. There 
is an $n>0$ such that $K \subseteq  C_n$. Let $L=C_{n+1}$. Then the 
incompressibility of $F_n$ in $M-Int \, C_0$ implies that any 
null-homotopy in $M-J$ of a loop in $M-L$ can be cut off on $F_n$ 
so as to obtain a null-homotopy in $M-K$. 

The $\bold P^2$-irreducibility of $M - \Inte C_0$ implies that $M$ is 
eventually $\bold P^2$-\irr. 

To show that $M$ is not almost compact it suffices to show that the 
fundamental group of some component of $M-C_0$ is not finitely generated. 
In fact, this is the case for every component $V$ of $M-C_0$. We may 
assume that $H_1(V)$ is finitely generated. Then for all sufficiently 
large $n$ one has that the intersection of each component of 
$V \cap X_{n+1}$  with each component of $V \cap X_{n+2}$ is connected. 
Since no component of $X_{n+1}$ has the form $F \times [0,1]$ 
described above it follows from Theorem 10.2 of \cite{10} that for every component 
$Y$ of $X_{n+1}$ and component  $S$ of $Fr \, Y$ we have that 
$\pi_1(S) \rightarrow \pi_1(Y)$ is a non-surjective monomorphism. 
It then follows  that $\pi_1(V)$ is an infinite non-trivial free product with 
amalgamation, hence is not finitely generated. 
\qed \enddemo   

Now suppose $M$ is a connected \3m\ with a finite number 
$\mu > 0$ of ends whose boundary consists of a finite number $\nu \geq 0$ 
of disjoint planes $E^i$. An exhaustion $C$ for $M$ is {\bf nice} 
if $C_n \cap \bd M$ consists of a single disk in each $E^i$, 
$X_{n+1}$ is $\bold P^2$-\irr, \birr, and \ann, each component of $F_n$ has 
negative Euler characteristic, each orientable component of $F_n$ has 
positive genus, and $M-Int \, C_n$ has $\mu$ components for all 
$n \geq 0$. Note again that $M$ need not be $\bold P^2$-\irr\ or \irr. 
            
\proclaim{Lemma 1.3} Let $C$ be a nice exhaustion for the connected, 
non-compact 3-man\-i\-fold $M$.  Then the following conditions hold. 
\roster
\item $C$ is a good exhaustion for $M$. 
\item $X_{n+1}\cap \bd M$ consists of one annulus in each component 
of $\bd M$ and is \inc\ in $X_{n+1}$. 
\item $M-C_0$, $M-Int\,C_0$, and, for $n\geq 1$, $C_n-Int\,C_0$ 
are $\bold P^2$-\irr, \birr, and \ann.  
\item In $M-Int \, C_0$ one has that $\bd M-Int(C_0 \cap \bd M)$ is 
\inc\ and $F_n$ is \binc. 
\item Every subexhaustion of $C$ is nice. 
\item $M$ is \ann\ at infinity. 
\endroster
\endproclaim 

\demo{Proof} Clearly $\{C_n \cap E^i\}$ is an exhaustion for $E^i$ 
by concentric disks, and so $X_{n+1} \cap E^i$ is an annulus with 
one boundary component in $\bd F_n$ and the other in $\bd F_{n+1}$. 
Suppose $D$ is a compressing disk for $F_n \cup F_{n+1}$ in $X_{n+1}$. 
Since $X_{n+1}$ is \birr\ there is a disk $D\p$ in $\bd X_{n+1}$ with 
$\bd D\p=\bd D$ such that $D\p$ must contain some annulus $X_{n+1} \cap 
E_i$ and hence some component of $F_n \cup F_{n+1}$, contradicting the fact 
that no component of $\FF0$ is a planar surface. Thus $F_n \cup F_{n+1}$ 
is \inc\ in $X_{n+1}$. If some component of $X_{n+1}$ is a product 
$F \times \I$ with $F \times \{0\}$ in $F_n$ and $F \times \{1\}$ in 
$F_{n+1}$, then since $F$ has negative Euler characteristic there is 
an \inc\ product annulus in $F \times \I$ which is not \bparallel\ , 
contradicting the fact that $X_{n+1}$ is \ann. Thus $C$ satisfies 
(1) and (2). 

Now suppose $D$ is a compressing disk for the boundary of 
$Y_n=C_n-Int \, C_0$, where $n>1$. Let $F=F_1 \cup \cdots \cup F_{n-1}$. 
If $D \cap F =\emptyset$, then $D$ lies in $X_1$ or $X_n$, say $X_1$.  
Isotop $D$ so that $\bd D$ lies in $F_0$. Then $\bd D=\bd D\p$ for 
a disk $D\p$ in $\bd X_1$.  If $D\p$ does not 
lie in $\bd Y_n$, then it must contain a component of $F_1$,  
contradicting the positive genus condition. Thus we may assume that 
 $D \cap F \neq \emptyset$. Suppose $D$ does not meet both $F_0$ and 
$F_n$; say it misses $F_n$. Then $D$ can be isotoped so that 
$\bd D$ lies in $F_0$. By the incompressibility of the $F_i$ and the 
irreducibility of the $X_j$ we may, if necessary, apply a finite 
sequence of disk pushes to $D$ so that $D \cap F=\emptyset$. Thus 
we may assume that $D$ meets both $F_0$ and $F_n$. By a disk push 
argument similar to that just given one may assume that $D \cap F$ 
contains no simple closed curves. Assume further that $D$ has been 
isotoped so that among all such $D$ the number of components of 
$D \cap F$ is minimal. Now 
$D \cap F$ splits $D$ into partial disks. Among these is a halfdisk 
$H$  which meets $\bd D$ in $\bd_0 H$ and $F$ in $\bd_1 H$. 
We may 
assume $H \subseteq X_1$, and so $\bd_1 H \subseteq F_1$. Then $\bd H =
\bd H\p$ for a disk $H\p$ in $\bd X_1$. This implies that $\bd (\bd_1 H)$ 
lies in a single component $J$ of $\bd F_1$, and so $J \cap H\p$ is 
an arc which splits $H\p$ into a disk $H_1\p$ in $F_1$ and a disk $H_0\p$ 
in $\bd X_1-int\, F_1$. Since $X_1$ is \irr\ $H \cup H\p$ bounds a ball $B$ 
in $X_1$, and so there is a halfdisk push of $H$ across $B$ 
past $H_1\p$ which removes at least $\bd_1 H$ from the 
intersection, thereby contradicting minimality. Thus $Y_n$ is \birr. 
The $\bold P^2$-irreducibility of $Y_n$ and $M-Int\,C_0$ follows from that of the 
$X_j$ together with the incompressibility of the $F_i$. This implies 
that $M-C_0$ is $\bold P^2$-\irr\ as well.

Suppose $D$ is a compressing disk for $\bd (M-Int\,C_0)$.  
Then for some $n$ one has $D\subseteq (Y_n-F_n)$, and so 
$\bd D= \bd D\p$ for a disk $D\p$ in $\bd Y_n$. 
If $D\p$ does not lie in $\bd(M-Int\,C_0)$, then it must 
contain a component of $F_n$. But this contradicts the 
positive genus condition, and so $M-Int\,C_0$ is \birr. 
The positive genus condition applied to $F_0$ now implies that 
$\bd M - Int(C_0 \cap \bd M)$ is \inc\ in $M-Int\,C_0$ and hence 
in $M-C_0$. 

Suppose $D$ is a $\bd$-compressing halfdisk in $M-Int \, C_0$ 
for some $F_n$. 
So $D \cap F_n=\bd_1 D$ and $D \cap \bd(M-Int \, C_0)=\bd_0 D$. 
Now either $D \subseteq Y_n$ or $D \cap Y_n
=\bd_1 D$. In the first case the $\bd$-irreducibility of $Y_n$ implies that 
$\bd D = \bd D\p$ for a disk $D\p$ in $\bd Y_n$. It follows that 
$\bd (\bd_1 D)$ 
lies in a component $J$ of $\bd F_n$ and that $J \cap D\p$ is an arc 
which splits $D\p$ into a disk in $F_n$ and a disk in 
$\bd Y_n -int \, F_n$; thus $\bd_1 D$ is \bparallel\ in $F_n$. In the second 
case $D$ lies in some $Y\p=C_{n+k}-Int \, C_n$, and one applies the 
same argument to this manifold. 

We have now established (4) and the portions of (3) not dealing 
with annuli. 
So suppose $A$ is an incompressible proper annulus in $Y_n$. 
Assume that the number of components of $A \cap F$ is minimal. 
Then a disk push argument shows that it contains no simple closed 
curve components which bound disks on $A$. 
Suppose $\al$ is a component of $A \cap F$ which is an arc. 
If $\al$ is \bparallel\ in $A$, then we may assume that $\al$ is 
innermost on $A$, hence $\al=\bd_1 D$ for a proper halfdisk $D$ on $A$ 
such that $D \cap F =\al$. Since $F$ is \binc\ 
in $M-Int \, C_0$ one has that $\al$ is \bparallel\ in $F$, 
and so $\al=\bd_1 D\p$ for a proper halfdisk 
$D\p$ in $F$. Then the disk $D \cup D\p$ is \bparallel\ in $Y_n$ 
via a $3$-ball $B$. 
Thus there is a halfdisk push of $D$ across $B$ past $D\p$ which 
 removes at least $\al$ from $A \cap F$, 
contradicting minimality. Thus $\al$ is a spanning arc in $A$. Since each 
$F_j$ separates $Y_n$ there must be two such arcs $\al$ and $\al\p$ such 
that $\al \cup \al\p=\bd_1 D$ for a proper band $D$ in $A$ such that 
$D \cap F = \al \cup \al\p$. 
Now $D$ lies in some $X_j$, and so $\bd D=\bd D\p$ for a disk $D\p$ in 
$\bd X_j$. Now $D \cup D\p$ bounds a 3-ball $B$ in $X_j$. 
Then $D\p \cap F$ consists of one or two disks. If it is a 
single disk $D\pp$, then $D\p$ is the union of $D\pp$ and two disks 
in $\bd Y_n$. There is an isotopy (a {\bf band push}) which moves 
the disk $D$ across $B$ past $D\pp$ which removes 
at least $\al \cup \al\p$ from $A \cap F$, contradicting minimality. 
If $D\p \cap F$ consists of two  disks then there are $\bd$-compressing 
halfdisks for $A$ in $Y_n$,  
and hence $A$ is \bparallel\ in $Y_n$. (See Lemma 2.2 of \cite{14}.) Thus 
we may assume $A \cap F$ contains no spanning arcs of $A$. 
Suppose $J$ is a simple closed curve component of $A \cap F$ which is 
non-contractible in $A$. Then $A$ is \bparallel\ in $A$ via an annulus 
$A_0$ in $A$; we may assume that $A_0 \cap F = J$. Then $A_0$ is a 
proper \inc\ annulus in some $X_j$. Since $X_j$ is \ann\ $A_0$ is 
parallel in $X_j$ to an annulus $A_1$ in $\bd X_j$ with $\bd A_1 = 
\bd A_0$. The parallelism $T$ is a solid torus with $\bd T=A_0 \cup A_1$. 
There is a component $A_2$ of $A_1 \cap F$ which is an annulus one of 
whose boundary components is $J$. Then there is an isotopy supported 
in a regular neighborhood of $T$ which moves $A_0$ across $T$ and past 
$A_2$, thereby removing at least $J$ from $A \cap F$ and thereby 
contradicting minimality. Thus we assume that $A \cap F = \emptyset$. 
Therefore $A$ lies in some $X_j$ and is 
parallel in $X_j$ to an annulus $A\p$ in $\bd X_j$ whose boundary misses 
that of $F$. By the positive genus condition $A\p \cap F=\emptyset$, and 
so $A\p \subseteq \bd Y_n$. Thus $Y_n$ is \ann. 

Now suppose $A$ is an \inc\ proper annulus in $M-Int\,C_0$. Then 
for some $n$ one has $A \subseteq (Y_n-F_n)$, and so $A$ is parallel 
to an annulus $A\p$ in $\bd Y_n$. If $A\p$ does not lie in 
$\bd(M-Int\,C_0)$, then it must contain a component of $F_n$. 
But this contradicts the positive genus condition, and so $M-Int\,C_0$ 
is \ann. The positive genus condition applied to $F_0$ implies that 
$M-C_0$ is \ann. This completes the proof of (3). 

Clearly the properties we have proven for $Y_n$ also hold for 
$C_n - Int \, C_m$ whenever $n>m$. It follows from this that every 
subexhaustion of $C$ is nice, thus establishing (5). 

Finally we note that given any compact subset $J$ of $M$ we can choose 
an $m > 0$ such that $J \subseteq C_m$. Then by (5) and (3) we have 
that $M-\Inte C_m$ is \ann. Suppose $A$ is an \inc\ proper annulus in 
$M-C_m$. Then $A$ is parallel in $M-\Inte C_m$ to an annulus $A\p$ in 
$\bd(M-\Inte C_m)$. Since $F_m$ is \inc\ in $M-\Inte C_m$ and has no 
disk or annulus components $A\p$ must lie in $\bd(M-C_m)$. Thus $M$ is 
\ann\ at infinity. \qed \enddemo

\head 2. Removing Trivial Intersections  \endhead

Given two compact, proper, incompressible surfaces in an irreducible 
3-manifold it is a standard result that one can ambiently isotop one of 
the surfaces so that the two surfaces are in general position 
and no simple closed curve component of their intersection is 
trivial, i.e\. bounds a disk in one (and hence both) of the surfaces. 
This isotopy consists of a finite series of disk pushes. In the case 
of non-compact surfaces one is faced with the possibility of an infinite 
series of disk pushes, which might not converge to an ambient isotopy. 
This section gives, in a slightly more general setting, conditions 
under which the isotopy can be accomplished, 
treating as well the removal of \bparallel\  arcs in \binc\  surfaces. 

\proclaim{Proposition 2.1} Let $M$ be a connected, \irr, non-compact 
\3m\  which is not \homeo\  to $\RRR$. Let $\Cal P$ and $\Cal Q$ be 
proper surfaces in $M$ which are in general position. Let $\JJ$ be a union 
of simple closed curve components of $\PP \cap \QQ$. Assume that the 
following conditions are satisfied.
\roster
\item No component of $\Cal P$ or of $\Cal Q$ is a 2-sphere. 
\item Each component $J$ of $\Cal J$ bounds 
a disk $D(J)$ on $\Cal P$ and a disk $G(J)$ on $\Cal Q$. 
\item There is no infinite sequence $\{J_m\}$ of distinct components of 
$\Cal J$ such that either $D(J_m) \subseteq int \, D(J_{m+1})$ for 
all $m$ or $G(J_m) \subseteq int \, G(J_{m+1})$ for all $m$, i.e. 
there is no {\bf infinite nesting} on $\Cal P$ or on $\Cal Q$ among 
the components of $\Cal J$. 
\endroster 
Then there is an 
ambient isotopy of $\Cal P$ in $M$, fixed on $\bd M$, which takes $\Cal P$ 
to a surface $\Cal P\p$ such that $\Cal P\p$ and $\Cal Q$ are in general 
position and $(\Cal P\p \cap \Cal Q) \subseteq (\Cal P \cap \Cal Q)-\Cal J$. 
Moreover, the isotopy is fixed on $\Cal P\p \cap \Cal Q$. 
\endproclaim

\demo{Proof} Since neither $\Cal P$ nor $\Cal Q$ contains 2-spheres 
the disks $D(J)$ and $G(J)$ are unique. Call $D(J)$ {\bf maximal} if 
there is no $D(J\p)$ such that $D(J) \subseteq int \, D(J\p)$. These 
maximal disks are disjoint and their union contains all the $D(J)$. 
Let $\{D(J_i)\}$ be the set of all maximal disks. Let $D_i$ be a regular 
neighborhood of $D(J_i)$ in $\Cal P$ chosen so that $D_i \cap \Cal Q =
D(J_i) \cap \Cal Q$ and distinct $D_i$ are disjoint. Let $\Cal D =
\cup D_i$. Similar remarks apply to the maximal disks on $\Cal Q$, 
yielding $\Cal G =\cup G_j$. It suffices to ambiently isotop $\Cal P$ 
to $\Cal P\p$ such that, denoting the image of $\Cal D$ under the isotopy 
by $\Cal D\p$, we have that $\Cal P\p$ and $\Cal Q$ are in general position,  
$\Cal D\p \cap \Cal Q = \emptyset$, 
$((\Cal P\p - \Cal D\p)\cap \Cal Q)$ 
is a union of components of 
$((\Cal P - \Cal D)\cap \Cal Q)$, and the isotopy is fixed on 
$((\Cal P\p - \Cal D\p)\cap \Cal Q)$.

Choose an exhaustion $C$ for $M$ such that for each $n \geq 0$ and each 
$m > n$ one has that $C_n$ and $C_m - Int \, C_n$ are \irr\.  
This is possible because $M$ is \irr\  and is not \homeo\  to $\RRR$. 
We may assume that $\Cal G \cap \Cal F_0 =\emptyset$ because the 
$G_j$ are disjoint disks in the interior of $M$. By passing to a subexhaustion of $C$ 
we may also assume that $Int \, C_{n+1}$ contains all those $D_i$ 
which meet $F_n$. Let $Y_0=Int \, C_1$, and let 
$Y_n=(Int \, C_{n+1})-C_{n-1}$ for $n \geq 1$. Thus if $D_i$ meets 
$F_n$, then $D_i \subseteq Y_n$. 

Suppose $n$ is even, $D_i$ meets $F_n$, and $D_i$ meets $\Cal G$. 
There is a disk $D$ in $D_i$ such that $D \cap \Cal G=\bd D$. Note 
that $\bd D$ need not be a component of $\Cal J$ and that $int \, D$ 
need not be disjoint from $\Cal Q$. Let $D_+$ be a regular neighborhood 
of $D$ in $\Cal P$ chosen so that $D_+ \cap \Cal Q = D \cap \Cal Q$. 
Let $G$ be the disk on $\Cal G$ 
bounded by $\bd D$. For a subset $A$ of $M$ let $A^*$ denote the image 
of $A$ under the disk push of $D$ past $G$ across the 3-ball $B$ in 
$Y_n$ bounded by $D \cup G$. 
We may assume that $\bd D_+$ lies on $\bd B^+$, where $B^+$ is the 
regular neighborhood of $B$ supporting the isotopy. 
Then $\Cal P^*$ and $\Cal Q$ are in general position, 
$(\Cal D^*-D_+^*) \cap \Cal Q$ and 
$((\Cal P^*- \Cal D^*) \cap \Cal Q)$ are unions of components of, 
respectively, $(\Cal D-D_+) \cap \Cal Q$ and 
$((\Cal P -\Cal D)\cap \Cal Q)$,  and $D_+^* \cap \Cal Q = \emptyset$. 

Now suppose $K$ is a component of $\Cal P \cap \Cal Q$. 
If $K$ lies outside $B$, then we may assume that the isotopy is fixed 
on $K$, and so $K$ is a  component of $\Cal P^* \cap \Cal Q$. 
If $K$ meets $B$, then it is a component of $D \cap \Cal Q$ or of 
$\Cal P \cap G$, and so it is removed by the isotopy. Since all the 
components of $\Cal P^* \cap \Cal Q$ are components of $\Cal P \cap \Cal Q$, 
we have that the isotopy is fixed on $\Cal P^* \cap \Cal Q$. 

We next consider the effect of the isotopy on other disks $D_k$. 
Suppose $D_k \cap F_n = \emptyset$. If $D_k \cap B = \emptyset$, then 
we may assume that $D_k^*=D_k$. If $D_k \cap B \neq \emptyset$, then 
$D_k \cap G \neq \emptyset$. The isotopy might move $D_k \cap B$ across 
some portion of $F_n$, but since $(D_k \cap B)^*$ lies within a regular 
neighborhood of $G$ and $G \cap F_n = \emptyset$ we may assume that 
$D_k^* \cap F_n = \emptyset$. Thus the number of disks $D_k^*$ in 
$\Cal P^*$ which meet $F_n$ is no greater than the number of disks 
$D_k$ in $\Cal P$ which meet $F_n$. Therefore after performing a finite 
sequence of these isotopies we may assume that if $D_i \cap F_n  
\neq \emptyset$, then $D_i \cap \Cal G = \emptyset$. Since these isotopies 
are supported in $Y_n$ we may do this simultaneously for all even $n$. 

Now for $n$ odd we still have that if $D_i \cap F_n \neq \emptyset$, 
then $D_i \subseteq Y_n$. Thus if $D_i$  is a component of $\Cal D$ such 
that $D_i \cap \Cal G \neq \emptyset$, then $D_i \subseteq Y_n$ 
for some odd $n$. By performing a finite sequence of disk pushes in each 
such $Y_n$ we obtain the desired surface $\Cal P\p$. \qed \enddemo

\proclaim{Corollary 2.2} Let $M$ be a connected, \irr, non-compact \3m\  
which is not \homeo\  to $\RRR$. Let $\Cal P$ and $\Cal Q$ be proper,  
\inc\  surfaces in $M$ such that no component of $\Cal P$ or of $\Cal Q$ 
is a 2-sphere or a trivial plane. Assume that there do not exist plane components 
$P$ of $\Cal P$ and $Q$ of $\Cal Q$ on both of which there is infinite nesting among 
the components of  $P \cap Q$.  Then there is an ambient isotopy of 
$\Cal P$ in $M$, fixed on $\bd M$, which takes $\Cal P$ to a surface 
$\Cal P\p$ such that $\Cal P\p$ and $\Cal Q$ are in general position 
and no simple closed curve component of $\Cal P\p \cap \Cal Q$ bounds 
a disk on $\Cal P\p$ or on $\Cal Q$. This isotopy is fixed on $\Cal P\p \cap 
\Cal Q$. \endproclaim  

\demo{Proof} Let $\Cal J$ be the set of all simple closed curve 
components of $\Cal P \cap \Cal Q$ which bound disks on $\Cal P$ 
(or equivalently on $\Cal Q$.) Infinite nesting on $\Cal P$ or 
on $\Cal Q$ among the elements of $\Cal J$ implies by Lemma 1.1 (1) 
that either one of these surfaces has a component which is a trivial plane or 
there are two plane components $P$ and $Q$ as above. 
\qed \enddemo

We now consider the removal of trivial arcs of intersection. 

\proclaim{Proposition 2.3} Let $M$ be a connected, irreducible, 
non-compact \3m\  
which has non-empty boundary and is not \homeo\  to $\halfspace$. 
Let $\Cal P$ and $\Cal Q$ be proper surfaces in $M$ which are in general 
position. Let $\AAA$ be a union of components of $\PP \cap \QQ$, each of 
which is an arc. Let $R$ be an end-proper surface in $\bd M$. Assume that the 
following conditions are satisfied. 
\roster 
\item No component of $\Cal P$ or of $\Cal Q$ is a 2-sphere or a disk. 
\item Each component $\al$ of $\Cal A$ is \bparallel\  
in $\Cal P$ across a halfdisk $D(\al)$ and is \bparallel\  in $\Cal Q$ across 
a halfdisk $G(\al)$. 
\item There is no infinite sequence $\{\al_m\}$ of 
distinct components of $\Cal A$ such that either 
$D(\al_m) \subseteq Int_{\Cal P}D(\al_{m+1})$ for all $m$ or 
$G(\al_m) \subseteq Int_{\Cal Q}G(\al_{m+1})$ for all $m$, 
i.e. there is no {\bf infinite nesting} on $\Cal P$ or on $\Cal Q$ among 
the components of $\Cal A$. 
\item $\bd \PP \cup \bd \QQ$ lies in $\inte R$.
\item $R$ is \inc\ in $M$. 
\item Each component $J$ of $\Cal J$ bounds a disk 
$G(J)$ in $\Cal Q$, where  $\Cal J$ is the union of all those simple 
closed curve components of $\Cal P \cap \Cal Q$ which lie in some 
$D(\al)$. 
\item  There is no infinite nesting on $\Cal Q$ 
among the components of $\Cal J$. 
\endroster
Then there is an ambient isotopy of $\Cal P$ in $M$, fixed on 
$(\bd M)-\inte R$ which 
takes $\Cal P$ to a surface $\Cal P\p$ such that $\Cal P\p$ and $\Cal Q$ 
are in general position and $(\Cal P\p \cap \Cal Q) \subseteq (\Cal P 
\cap \Cal Q) - (\Cal A \cup \Cal J)$. Moreover, the isotopy is fixed on 
$\Cal P\p \cap \Cal Q$. \endproclaim

\demo{Proof} Since neither $\Cal P$ nor $\Cal Q$ has a disk component 
the halfdisks $D(\al)$ and $G(\al)$ are unique. As before call a 
halfdisk $D(\al)$ {\bf maximal} if there is no $D(\al\p)$ such that 
$D(\al) \subseteq Int_{\Cal P}D(\al\p)$. By hypothesis each $D(\al)$ 
lies in some maximal halfdisk and these maximal halfdisks are disjoint. 
Let $\{D(\al_i)\}$ be the set of maximal halfdisks. Let $D_i$ be a 
regular neighborhood of $D(\al_i)$ in $\Cal P$ chosen so that 
$D_i \cap \Cal Q=D(\al_i) \cap \Cal Q$ and distinct $D_i$ are disjoint. 
Let $\Cal D=\cup D_i$. The maximal halfdisks in $\Cal Q$ yield 
$\Cal G=\cup G_j$. It suffices to ambiently isotop $\Cal P$ to $\Cal P\p$ 
such that, denoting  the image of $\Cal D$ under the isotopy by $\Cal D\p$, 
we have that  $\Cal P\p$ and $\Cal Q$ are in general position, 
$\Cal D\p \cap \Cal G = \emptyset$, $((\Cal P\p - \Cal D\p)\cap \Cal Q)$ is a 
union of components of  $((\Cal P - \Cal D)\cap \Cal Q)$, and the isotopy 
is fixed on $((\Cal P\p - \Cal D\p)\cap \Cal Q)$.

Since $\Cal J \subseteq \Cal D$ there is no infinite nesting on 
$\Cal P$ among the components of $\Cal J$.  
Thus by Proposition 2.1 we may assume $\Cal J=\emptyset$. 

As before we can choose an exhaustion  $C$ for $M$ such that for 
each $n\geq 0$ and each $m>n$ one has that $C_n$ and $C_m-Int \, C_n$ 
are irreducible. Since the $G_j$ are disjoint disks each meeting $\bd M$ in 
a single arc we may assume that $\Cal G \cap \Cal F_0 = \emptyset$. 
We may also assume that $\bd \FF0$ and 
$\bd R$ are in general position, and that if $D_i$ meets $F_n$, then 
$D_i \subseteq Int \, C_{n+1}$. Note that since $M$ is 
not \homeo\  to $\halfspace$  we may assume that no $C_n$ lies in a 
3-ball 
in $M$ which meets $\bd M$ in a single disk. We claim that $C_n \cap R$ 
and $(C_m-Int \, C_n) \cap R$ are \inc\  in $C_n$ and $(C_m-Int \, C_n)$ 
respectively. If $D$ is a compressing disk for $C_n \cap R$ in $C_n$, 
then $\bd D=\bd D\p$ for a disk $D\p$ in $R$. If $D\p$ does not lie in 
$C_n \cap R$, then it must meet $M-C_n$. Since 
$M$ is \irr\  $D \cup D\p$ bounds a 3-ball containing a component of 
$M-C_n$, contradicting the fact that $M-C_n$ has no bounded components. 
If $D$ is a compressing disk for $(C_m-Int \, C_n) \cap R$ in 
$C_m-Int \, C_n$, then $\bd D = \bd D\p$ for a disk $D\p$ in $R$ which 
must meet $C_n$ or $M-C_m$. Let $B$ be the 3-ball in $M$ bounded by 
$D \cup D\p$. Then $D\p$ does not meet $C_n$, for otherwise $C_n$ would 
lie in $B$. So $D\p$ meets $M-C_m$ and hence a component of $M-C_m$ 
is contained in $B$, contradicting the fact that $M-C_m$ has no 
bounded components. 

We now proceed as in the proof of Proposition 2.1. We let $Y_0=Int \, C_1$ 
and $Y_n=(Int \, C_{n+1})-C_{n-1}$ for $n \geq 1$. Let $R_n=Y_n \cap R$. 
Suppose $D_i$ meets both $F_n$ and $\Cal G$. The role of an innermost disk 
is now given to 
an innermost halfdisk, i.e. a proper halfdisk $D$ in $\Cal P$ which lies 
in $D_i$ such that $\bd_1 D = D_i \cap \Cal G$, although $\bd_1 D$ need 
not be a component of $\Cal A$ and $Int_{\Cal P}D$ need not be disjoint from 
$\Cal Q$. There is a unique proper halfdisk $G$ in $\Cal Q$ such that $G$ 
lies in $\Cal G$ and $\bd_1 G = \bd_1 D$. Then $D \cup G$ is a proper 
disk in $Y_n$ with $\bd(D \cup G)$ in $R_n$. Since $R_n$ is \inc\  in 
$Y_n$ there is a disk $D\p$ in $R_n$ with $\bd D\p=\bd(D \cup G)$. Since 
$Y_n$ is \irr\  $D \cup G \cup D\p$ bounds a 3-ball $B$ in $Y_n$. 
A halfdisk push of $D$ across $B$ past $G$ removes at least $\bd_1 D$ 
from the intersection, adds no new components of intersection, either fixes or 
removes each component of $\Cal P \cap \Cal Q$, and does not increase the 
number of $D_k$ meeting $F_n$. As before we remove all intersections of 
$D_i$ with $\Cal G$ for those $D_i$ which meet $F_n$ for $n$ even and 
then remove all remaining intersections of $\Cal D$ with $\Cal G$. 
\qed \enddemo

\proclaim{Corollary 2.4} Let $M$ be a connected, irreducible, 
non-compact \3m\  
which has non-empty boundary and is not homeomorphic to $\halfspace$. 
Let $\Cal P$ and $\Cal Q$ be proper, incompressible, \binc\  surfaces 
in $M$ 
such that no component of $\Cal P$ or of $\Cal Q$ is a 2-sphere, a disk,  
a trivial plane, or a trivial halfplane. Assume that there do not exist components 
$P$ of $\Cal P$ and $Q$ of $\Cal Q$ which are either both planes or both halfplanes 
on both of which there is infinite nesting among the components of $P \cap Q$. 
Suppose $R$ is an end-proper surface 
in $\bd M$ which is \inc\ in $M$ and whose interior contains 
$\bd \Cal P \cup \bd \Cal Q$.  
Then there is an ambient isotopy of $\Cal P$ in $M$, fixed on  
$(\bd M)-\inte R$, which 
takes $\Cal P$ to a surface $\Cal P\p$ such that $\Cal P\p$ and $\Cal Q$ 
are in general position, no simple closed curve component of $\Cal P\p 
\cap \Cal Q$ bounds a disk on $\Cal P\p$ or on $\Cal Q$, and no  
component of $\Cal P\p \cap \Cal Q$ is an arc which is \bparallel\ 
in $\Cal P\p$ or in $\Cal Q$. This isotopy is fixed on $\Cal P\p \cap \Cal Q$. 
\endproclaim

\demo{Proof} First apply Corollary 2.2 to remove all trivial simple 
closed curve components of $\Cal P \cap \Cal Q$. Then let $\Cal A$ 
be the set of all those components of $\Cal P \cap \Cal Q$ which are 
\bparallel\ arcs in $\Cal P$ (or equivalently in $\Cal Q$.) Infinite 
nesting among the components of $\Cal A$ implies by Lemma 1.1 (2) 
that either one of these surfaces has a component which is a trivial 
halfplane or there are halfplane components $P$ and $Q$ as above. \qed \enddemo

\head 3. Aplanarity and total acylindricality  \endhead

The goal of this section is to show that a connected, non-compact  
3-manifold which possesses a nice exhaustion must be totally acylindrical 
and, if it is \irr, must be aplanar. These results are basically due to 
Winters and are contained, either explicitly or implicitly, in his thesis 
[19], where they sometimes 
appear in a slightly different form and context. We include proofs of them 
here for several reasons. First, the paper [20] containing the relevant portions 
of [19] has not yet been published, and so giving proofs here will make the 
argument of the present paper more complete. Second, the proof given here of 
Lemmas 3.1 and 3.2 is somewhat different from that of the corresponding 
Lemma X.1 of [19] in that it applies the general machinery for removing 
trivial intersections developed in the previous section of this paper rather 
than the direct arguments of [19]; this remark also applies to the proof of 
Theorem 3.5. Third, Lemma 3.3 and Theorems 3.4 and 3.5 are not 
stated explicitly in [19] in the forms we shall need, although their proofs either 
are contained in or can be easily deduced from the proofs of Lemmas II.2, II.3, 
and XII.3 of [19]. Finally, the terminology and the organization of the proof 
of aplanarity given in 
this section establish the background and conceptual framework for the 
more difficult analysis of partial planes in the next section.

In the first three lemmas 
$M$ is a connected non-compact 3-manifold and 
$C$ is an exhaustion for $M$. 
Whenever $J$ is a simple closed curve on a plane $P$, let $D(J)$ be the 
disk on $P$ bounded by $J$. Let $n_0 > 0$. A proper plane $P$ is in 
{\bf $n_0$-standard position} with respect to $C$ if $P$ is in general 
position with respect to $\FF0$, $P \cap \FFn0$ is a sequence 
$\{J_m\}$ of simple closed curves such that 
$\{D(J_m)\}$ is an exhaustion for $P$, 
$D(J_0)$ is a proper disk in $\Cn0$, and 
$(P \cap C_0) \subseteq int \, D(J_0)$. 
Note that if this is the case and $n_1 > n_0$, then $P$ is also in 
$n_1$-standard position with respect to $C$. 

\proclaim{Lemma 3.1} Suppose $C$ is a good exhaustion for $M$ and   
$P$ is a proper plane in $M$. Then for some $n_0 > 0$ 
one has that $P$  is ambient isotopic to a plane which is in 
$n_0$-standard position with respect to $C$. 
\endproclaim

\demo{Proof} First ambiently isotop $P$ so that it is in general position 
with respect to $\Cal{F}_0$. Since $P$ is proper there is a disk $D$ in $P$ 
which contains $P \cap C_0$ in its interior. (If $P \cap C_0 = \emptyset$, 
let $D$ be any disk in $P$.) Since $C$ is an exhaustion and $D$ is compact 
there is an $n_0 > 0$ such that $D \subseteq Int \, \Cn0$. There is then a 
component $J$ of $P \cap \Fn0$ such that $D \subseteq int \, D(J)$. 
Let $J_0$ be the innermost such component of $P \cap \Fn0$, i.e. there 
is no component $J$ of $P \cap \Fn0$ such that $D \subseteq int \, D(J) 
\subseteq int \, D(J_0)$.

Now let $J$ be a component of $D(J_0) \cap \Fn0$ other than $J_0$. 
Assume that $J$ is innermost on $P$ among such curves. Then 
$D(J) \cap \Fn0=J$. Moreover $D(J) \cap D = \emptyset$, and  so 
$D(J)$ lies in $M - Int \, C_0$. Since $\Fn0$ is incompressible in 
$M - Int \, C_0$ there is a disk $D\p$ in $\Fn0$  bounded by $J$. 
Since $M - Int \, C_0$ is irreducible the  2-sphere $D(J) \cup D\p$ bounds 
a 3-ball $B$ in $M - Int \, C_0$, and so one can perform a disk push 
of $D(J)$ across $B$ past $D\p$. This ambient isotopy of $P$ is 
supported in $M - Int \, C_0$,  
removes $J$ from the intersection, and adds no new components. 
We continue performing such isotopies until $D(J_0)$ becomes a proper 
disk in $\Cn0$. 

Suppose $J$ is a component of $P \cap \Fn0$ such that $D(J)$ does 
not contain $D(J_0)$. Then $D(J)$ lies in $M-Int \, C_0$ and so  
a sequence of disk pushes similar to that described above removes 
all such $J$ from the intersection. 

Let $\Yn0= M-Int \, \Cn0$ and $\Cal P=P \cap \Yn0$. 
Then $\Yn0$ is \irr\ and \birr\, and the set of components of $\Cal P$ 
consists of a half-cylinder (homeomorphic to $S^1 \times [0,\infty)$) 
and possibly a finite collection of annuli, all of which are proper in 
$\Yn0$. Let $\Cal Q={\Cal F}_{n_0+1}$. Let $\Cal J$ be the union of 
all those components $J$ of $\Cal P \cap \Cal Q$ such that $D(J)$ does 
not contain $D(J_0)$. Then $\Yn0$, $\Cal P$, $\Cal Q$, and $\Cal J$ 
satisfy the hypotheses of Proposition 2.1, and so there is an ambient 
isotopy in $\Yn0$, fixed on $\Fn0$, which removes $\Cal J$ from the 
intersection and adds no new components to it. Since the isotopy is 
fixed on those intersection curves which are not removed they all bound 
disks which contain $D(J_0)$. The isotopy thus extends by 
the identity isotopy on $\Cn0$ to an ambient isotopy in $M$ which 
carries $P$ to a plane $P\p$ which is in $n_0$-standard position with 
respect to $C$. \qed \enddemo 

A proper plane $P$ is in {\bf non-trivial $n_0$-standard position} 
with respect to an exhaustion $C$ if it is in $n_0$-standard position 
with respect to $C$ and no component of $P \cap \FFn0$ bounds a disk 
in $\FFn0$. 

\proclaim{Lemma 3.2} Let $C$ be a good exhaustion for $M$, and let $P$ be 
a non-trivial proper plane in $M$. 
Suppose $M$ is irreducible. Then for some $n_0 > 0$ 
one has that $P$ is ambient 
isotopic to a plane which is in non-trivial $n_0$-standard position with 
respect to $C$. \endproclaim 

\demo{Proof} By Lemma 3.1 we may assume that $P$ is in $k$-standard 
position with respect to $C$ for some $k > 0$. Then 
$P \cap \Cal{F}_k$ is a nested sequence $\{J_m\}$ of simple closed 
curves on $P$. If there is an increasing sequence $\{m_i\}$ such that 
$J_{m_i}$ bounds a disk $D\p(J_m)$ in $\Cal{F}_k$, then 
by Lemma 1.1 (1) we have that $P$ is trivial in $M$.  

Thus no such sequence exists, and so there is an $m_0 \geq 0$ such that 
for all $m \geq m_0$ $J_m$ does not bound a disk in $\Cal{F}_k$. 
There is then an $m_1 \geq m_0$ and an $n_0 \geq k$ such that 
$J_{m_1} \subseteq \Fn0$ and $D(J_{m_1}) \subseteq \Cn0$. Then $P$ is in 
non-trivial $n_0$-standard position with respect to $C$. 
\qed \enddemo 

A proper plane $P$ in $M$ is in {\bf $n_0$-monotone position} with 
respect to an exhaustion $C$ if it is in $n_0$-standard position 
with respect to $C$ and for $n \geq n_0$ one has that $P \cap X_{n+1}$ is 
an annulus with one boundary component in $F_n$ and the other in 
$F_{n+1}$. We apply the adjective ``non-trivial'' in the same sense 
as above. Again $P$ is in $n_1$-monotone position for all $n_1 > n_0$. 
We now assume that $M$ is a connected non-compact 3-manifold 
having finitely many ends and finitely many boundary components and that 
each boundary component is a plane. 

\proclaim{Lemma 3.3} Suppose $C$ is a nice exhaustion for $M$ and 
$P$ is a non-trivial proper plane in $M$. Assume that $M$ is 
irreducible. Then for some $n_0 > 0$ one has that $P$ is ambient 
isotopic to a plane which is in non-trivial $n_0$-monotone position 
with respect to $C$. \endproclaim 

\demo{Proof} By Lemma 3.2 we may assume that $P$ is in non-trivial 
$n_0$-standard position with respect to $C$ for some $n_0 > 0$. 
Then there 
is an exhaustion $\{D_m\}$ of $P$ with each $D_m$ a disk such that 
$(P \cap C_0) \subseteq int \, D_0$ and $D_0$ is a proper disk in $C_{n_0}$. 
Moreover, since $\bd D_m$ does not bound a disk in $\FFn0$ and $\FFn0$ is 
incompressible in $M-Int \, C_0$, one has that each annulus 
$D_{m+1}-int \, D_m$ 
is incompressible in the $X_{n+1}$ containing it. Denote $\bd D_m$ by 
$J_m$. 

Let $m_0=0$. Let $m_1$ be the smallest index for which $J_{m_1} \subseteq 
F_{n_0 +1}$ and $(P \cap \Cn0) \subseteq D_{m_1}$. Assume $m_1, \ldots, 
m_r$ have been defined. Let $m_{r+1}$ be the smallest index for which 
$J_{m_{r+1}} \subseteq F_{n_0 + r + 1}$ and $(P \cap C_{n_0 + r}) \subseteq 
D_{m_{r+1}}$. Then $\{D_{m_r}\}$ is an exhaustion for $P$. Let $A_{r+1}=
D_{m_{r+1}}-int \, D_{m_r}$. Suppose $(int \, A_{r+1}) \cap \FFn0 \neq \emptyset$. 
Call any component of this set a {\bf redundant intersection}. Since each 
$F_n$ separates $M$ the redundant 
intersections occur in pairs, with each component of a pair lying in the 
same component of $\FFn0$ and the pair forming the boundary of an annulus 
contained in $A_{r+1}$. Call such an annulus a {\bf redundant annulus}. 
There is then a redundant annulus $A$ which is innermost in the sense that 
its interior misses $\FFn0$. Then $A$ is a proper incompressible annulus in 
some $X_{n+1}$, and $\bd A \subseteq F_n$ or $\bd A \subseteq F_{n+1}$. 

Assume $\bd A \subseteq F_n$. Then $A$ is parallel in $X_{n+1}$ to an annulus 
$A\p$ in $\bd X_{n+1}$. Suppose $A\p$ does not lie in $F_n$. By the 
incompressibility of $F_n$, $F_{n+1}$, and $X_{n+1} \cap \bd M$ in $X_{n+1}$ 
the components of $A\p \cap F_n$, $A\p \cap F_{n+1}$ and $A\p \cap \bd M$ 
must be annuli. $A\p \cap F_{n+1} = \emptyset$ since its components would be 
components of $F_{n+1}$, which, since $C$ is nice, has no annulus 
components. Since $F_n$ has no annulus components $A\p \cap F_n$ has exactly 
two components, each of which is a collar on a boundary component of $F_n$. 
Thus $A\p \cap \bd M$ is an annulus whose boundary lies entirely in $F_n$, 
again contradicting the fact that $C$ is nice. Therefore $A\p$ lies in 
$F_n$ and so there is an ambient isotopy supported in a regular 
neighborhood of $X_{n+1}$ which reduces the number of redundant intersections. 
A similar argument holds for $\bd A \subseteq F_{n+1}$. Note that this process 
may move some redundant annuli and remove others. But it introduces no new 
redundant annuli and leaves any remaining redundant annuli lying in the 
union of the same set of $X_{k+1}$'s. 

Thus any one redundant annulus can be removed by an ambient isotopy 
of $P$ supported in $M- Int \, C_0$. However, since there may be 
infinitely many redundant annuli one must avoid performing an infinite 
sequence of these isotopies which fails to converge to an ambient isotopy 
of $M$. We proceed as follows.

Choose $n_1 > n_0$ so that if $A$ is any redundant annulus with 
$\bd A \subseteq \Fn0$, then $A$ lies in $(Int \, C_{n_1})-C_0$, and if $H$ is 
any redundant annulus with $\bd H \subseteq F_{n_1}$, then $H$ lies in 
$M-\Cn0$. Now suppose $n_0, n_1, \ldots, n_i$, $i > 0$, have been chosen. 
Choose $n_{i+1} > n_i$ so that if $A$ is any redundant annulus with 
$A \subseteq F_{n_i}$, then $A$ lies in $(Int \, C_{n_{i+1}})- C_{n_{i-1}}$, and 
if $H$ is any redundant annulus with $\bd H \subseteq F_{n_{i+1}}$, then $H$ 
lies in $M-C_{n_i}$. 

By an ambient isotopy supported in $(Int \, C_{n_1})-C_0$ remove all redundant 
intersections with $\Fn0$. For each even $i > 0$ perform an ambient isotopy 
supported in $(Int \, C_{n_{i+1}})-C_{n_{i-1}}$ which removes all redundant 
intersections with $F_{n_i}$. Since they have disjoint compact supports these 
isotopies give a single ambient isotopy supported in $M-C_0$. Now for each 
odd $i \geq 1$ perform an ambient isotopy which removes all redundant 
intersections with $F_{n_i}$. Again this gives a single ambient isotopy 
supported in $M-C_0$. 

One now has no redundant intersections with the $F_{n_i}$. For each $i > 0$ 
perform an ambient isotopy supported in $(Int \, C_{n_i})-C_{n_{i-1}}$ which 
removes all redundant intersections with those $F_j$ with $n_{i-1} < j < 
n_i$. This completes the removal of all redundant intersections and puts $P$ 
in non-trivial $n_0$-monotone position with respect to $C$. \qed \enddemo

\proclaim{Theorem 3.4} Let $M$ be a connected, irreducible, non-compact 
\3m\ which has a nice exhaustion. Then $M$ is aplanar. 
\endproclaim

\demo{Proof} Suppose $P$ is a nontrivial proper plane in $M$, and 
$C$ is a nice exhaustion for $M$. 
By Lemma 3.3 we may assume that $P$ is in non-trivial $n_0$-monotone  
position with respect to $C$. For simplicity of notation we shall reindex 
so that $n_0=0$. Then $P \cap C_n$ is a disk $D_n$, $\{D_n\}$ is an 
exhaustion of $P$, and $P \cap X_{n+1} = D_{n+1} - int \, D_n$ is an annulus 
$A_{n+1}$ with one boundary component in $F_n$ and the other in $F_{n+1}$. 
Moreover $J_n=\bd D_n$ does not bound a disk in $F_n$. 

It follows that each $A_{n+1}$ is parallel in $X_{n+1}$ to an annulus 
$A_{n+1}\p$ in $\bd X_{n+1}$. Since ${\Cal F}_0$ has no disk components and 
is incompressible in $M-Int \, C_0$ one has that each component of 
$A_{n+1}\p \cap (F_n \cup F_{n+1})$ is an annulus. Since ${\Cal F}_0$ has 
no annulus components this intersection consists of collars $G_{n+1}$ and 
$H_{n+1}$, respectively, on $J_n$ and $J_{n+1}$ in $A_{n+1}\p$. Let 
$R_{n+1}=A_{n+1}\p \cap \bd M$. Then $A_{n+1}\p=G_{n+1} \cup R_{n+1} 
\cup H_{n+1}$. The parallelism between $A_{n+1}$ and $A_{n+1}\p$ defines an 
embedding of $S^1 \times \I \times \I$ in $X_{n+1}$, with 
$S^1 \times \I \times \{0\}=A_{n+1}$, $S^1 \times \I \times \{1\}=R_{n+1}$, 
$S^1 \times \{0\} \times I = G_{n+1}$, and $S^1 \times \{1\} \times I = 
H_{n+1}$. Consider the analogous situation in $X_{n+2}$. If $G_{n+2} \neq 
H_{n+1}$, then $G_{n+2} \cup H_{n+1}$ is a component of $F_{n+1}$, 
contradicting the fact that ${\Cal F}_0$ has no annulus components. 
Therefore $G_{n+2}=H_{n+1}$, all the $R_{n+1}$ lie in the same component 
$E_i$ of $\bd M$, and so one can fit together the embeddings of 
$S^1 \times \I \times \I$ to get an end proper embedding of 
$S^1 \times [1,\infty) \times \I$ in $M-Int \, C_0$ with 
$S^1 \times [1,\infty) \times \{0\}=P-int \, D_0$, 
$S^1 \times [1,\infty) \times \{1\} =E_i-(E_i \cap C_0)$, and 
$S^1 \times \{1\} \times \I=G_1$. Now $D_0 \cup G_1$ is a proper disk in $M$ 
whose boundary is that of $E_i \cap C_0$.  
Since $M$ is irreducible the union of these disks  
bounds a 3-ball in $M$ which can be used to extend the product structure to 
obtain a parallelism between $P$ and $E_i$. \qed \enddemo

Now recall that $M$ is totally acylindrical if it admits no proper 
incompressible embeddings of a cylinder $S^1 \times \R$. 

\proclaim{Theorem 3.5} Let $M$ be a connected, non-compact 
\3m\  which has a nice exhaustion. Then $M$ is totally 
acylindrical. \endproclaim

\demo{Proof} Suppose $S=S^1 \times \R$ is a proper, \inc\ 
cylinder in $M$. Let $C$ be a nice exhaustion for $M$. Put $S$ in 
general position with respect to $\FF0$. Since $S$ is proper 
there exist $a$, $b \in \R$, $a<b$, such that the annulus 
$A=S^1 \times [a,b]$ contains $S \cap C_0$. (If $S \cap C_0=
\emptyset$, then choose $a<b$ arbitrarily.) There is then an 
$n_0 > 0$ such that $A \subseteq Int \, \Cn0$. There are 
components $J_0^+$ and $J_0^-$ of $S \cap \Fn0$, neither 
of which bounds a disk on $S$, such that $A$ is contained 
in the interior of the annulus $A_0$ on $S$ bounded by 
$J_0^+ \cup J_0^-$. We may assume that $A_0$ is an innermost 
such annulus, i.e. there are no components $J^+$, $J^-$ of 
$S \cap \Fn0$ which are non-contractible on $S$ such that the 
annulus bounded by $J^+ \cup J^-$ contains $A$ and is 
contained in, but does not equal, the annulus $A_0$. 

We now proceed as in the proof of Lemma 3.1. All the components 
of $A_0 \cap \Fn0$ other than $J^+$ and $J^-$ bound disks 
in $A_0$ which miss $A$ and lie in $M-Int\,C_0$. We 
perform a finite sequence of disk pushes which remove these 
curves. We then isotop the two half-cylinders composing 
$S-Int\,A_0$ to remove all the other components of 
$S \cap \Fn0$ which bound disks on $S$. 

Let $\Yn0=M-Int\,\Cn0$ and $\PP=S \cap \Yn0$. Then $\PP$ 
consists of two half-cylinders and possibly a finite collection of 
annuli, all of which are proper in $\Yn0$. Let $\QQ=\Cal F_{n_0+1}$. 
Let $\JJ$ be the union of all those components of $\PP \cap \QQ$ 
which bound disks on $S$ and hence on $\PP$. Since $Y_{n_0}$ is \irr\ 
Proposition 2.1 gives an ambient isotopy of $S$ in $M$ which 
removes $\JJ$ from the intersection, adds no new components, and 
is fixed on those which are not removed.  

Now each $X \cap X_{n+1}$ is a union of disjoint annuli. 
Since $S$ is proper this intersection is non-empty for all 
sufficiently large $n$ and in fact must contain an annulus 
$A\p$ running from $F_n$ to $F_{n+1}$. Since $X_{n+1}$ is 
\ann\ $A\p$ is parallel to an annulus $A\pp$ in $\bd X_{n+1}$. 
Now $A\pp$ must contain the annulus $E \cap X_{n+1}$ for 
some component $E$ of $\bd M$. Since $E$ is a plane it 
follows that $S$ is compressible in $M$, a contradiction. 
\qed \enddemo

\head 4. Strong Aplanarity: The Special Case  \endhead

Recall that a partial plane $P$ is a non-compact, simply-connected 
2-manifold with non-empty boundary.  In this section we show that 
given any proper partial plane $P$ in an irreducible 3-manifold $M$ 
which has a nice exhaustion there exists a collar on $\bd M$ which 
contains $P$. This is a strong condition on $M$ since such a 
$P$ cannot meet distinct components of $\bd M$. One should note, 
however,  that $P$ need not be \bparallel. In the next section this 
result will be extended to proper surfaces each of whose components 
is a partial plane. 

The general line of argument will be similar to that of the previous 
section. Instead of working with simple closed curves and the disks 
they bound we shall work with finite collections of proper arcs which 
cut off disks on $P$. Specifically, let $\al$ be the union of finitely 
many disjoint proper arcs $\al^1, \, \ldots \, , \, \al^k$ in a partial 
plane $P$. If for some proper partial disk $D$ in $P$ one has 
$\al=\bd_1 D$, then $\al$ is called a {\bf bounding arc system} in $P$, 
and $D$ is denoted by $D(\al)$. 

A proper partial plane $P$ is in {\bf \nzsp} with respect to 
an exhaustion $C$ 
if $P$ is in general position with respect to $\Cal F_0$, each component 
of $P \cap \FFn0$ is an arc, there is a sequence $\{D_m\}$ of 
well embedded partial disks in $P$ which is an exhaustion for  
$P$, $(P \cap C_0) \subseteq Int_P  D_0$, $D_0$ is a proper disk in 
$\Cn0$, and $P \cap \FFn0 = \cup_{m \geq 0} Fr_P  D_m$. 
In this case it is in $n_1$-standard position for all $n_1 > n_0$. 

\proclaim{Lemma 4.1} Let $C$ be a nice exhaustion for $M$, and let 
$P$ be a proper partial plane in $M$. Then for some $n_0 > 0$ 
one has that $P$ is ambient isotopic to a partial plane which is 
in \nzsp\ with respect to $C$. \endproclaim

\demo{Proof} First ambiently isotop $P$ so that it is in general position 
with respect to $\Cal{F}_0$. Since $\bd P \neq \emptyset$ and $P$ is proper 
$\bd M \neq \emptyset$. Since $P$ is proper there is a disk $D$ in $P$ 
such that $(P \cap C_0) \subseteq Int_P D$ and $D \cap \bd P \neq 
\emptyset$. (If $P \cap C_0 = \emptyset$, let $D$ be any disk in $P$ with 
$D \cap \bd P \neq \emptyset$.) Since $C$ is an exhaustion and $D$ is compact 
there is an $n_0 > 0$ such that $D \subseteq Int \, \Cn0$. There is a unique component 
$L$ of $P \cap C_{n_0}$ such that $D \subseteq Int_P L$. Then $L$ is a planar surface 
whose boundary consists of one simple closed curve $K$ which meets $\bd P$ and 
possibly some other simple closed curves $K_j$ which do not meet $\bd P$. Let 
$\be=K \cap \Fn0$. Then $\be$ is a bounding arc system in $P$ such that $D(\be)$ 
is the union of $L$ with the disks in $P$ bounded by the $K_j$.

If $D(\be)$ is not well embedded in $P$, then for some 
components $\be^i$ of $\be$ there exist proper halfdisk components $D^i$ 
of $P-Int_P D(\be)$ such that $\bd_1D^i=\be^i$.  Delete these $\be^i$ 
from $\be$ to obtain a new arc system $\al$. Then $\al$ is a bounding arc 
system such that  $D(\al)$ is well embedded in $P$, 
$\al \subseteq (P \cap \Fn0)$, and $D \subseteq Int_P D(\al)$. 
Note that $D(\al)$ need not lie in $C_{n_0}$. We shall next isotop $P$ so 
that afterwards $D(\al)$ does lie in $C_{n_0}$. To simplify the 
notation denote $D(\al)$ by $D_0$. 

Suppose $J$ is a simple closed curve component of $D_0 \cap \Fn0$. 
We may assume that $J$ is innermost on $P$ among such curves.  
Since $D \cap \bd P \neq \emptyset$ one has that 
$D(J) \cap D = \emptyset$, and so $D(J)$ lies in $M - Int \, C_0$. Since 
$M-Int \, C_0$ is irreducible and $\Fn0$ is incompressible in $M-Int \, C_0$, 
there is a disk push which removes $J$ from the intersection and adds 
no new components. Continue in this fashion until all such $J$ are removed. 

Now suppose $\ga$ is a component of $(Int \, D_0) \cap \Fn0$ which 
is a proper arc in $D_0$. We may assume that $\ga$ is innermost on $P$ 
among such arcs. Hence there is a proper halfdisk $H$ in $P$ with 
$\bd_1 H = \ga = H \cap F_{n_0}$.  We may further assume that $H$ 
does not contain $D$, and so $H$ lies in $M-Int \, C_0$. Since 
$F_{n_0}$ is \binc\ in $M-Int \, C_0$ there is a proper halfdisk 
$H\p$ in $F_{n_0}$ with $\bd_1 H\p = \ga$. Since 
$\bd M - Int(C_0 \cap \bd M)$ is \inc\ in $M-Int \, C_0$, the simple 
closed curve $\bd_0 H \cup \bd_0 H\p$ bounds a disk $H\pp$ in 
$\bd M - Int(C_0 \cap \bd M)$. Hence by the irreducibility of 
$M - \Inte C_0$ there is a halfdisk push of $H$ past $H\p$ across the 
ball bounded by $H \cup H\p \cup H\pp$ which removes 
$\ga$ from the intersection and adds no new components. We continue 
in this fashion until $D_0$ becomes a proper disk in $C_{n_0}$. 

We now consider components of $P \cap F_{n_0}$ which do not lie in 
$D_0$. As above we first remove all such components which are 
simple closed curves and then remove all those which are \bparallel\  
arcs in $P$. Note that afterwards $P \cap F_{n_0}$ splits $P$ into 
finitely many partial disks and partial planes. 

Let $Y_{n_0} = M - Int \, C_{n_0}$. We now have that $P \cap C_{n_0}$ 
consists of $D_0$ and possibly a finite number of other disks. 
$P \cap \Yn0$ consists of a finite number of partial planes and 
perhaps finitely many disks $D_j$, each of which meets $\Fn0$. Let 
$\PP$ be the union of the partial planes, and let 
$\QQ =\Cal F_{n_0+1}$. By an isotopy fixed on $\bd \Yn0$ we remove all 
simple closed curve intersections of the $D_j$ with $\QQ$. Then 
by an isotopy fixed on $\Fn0$ we remove all those components of 
$(\cup D_j) \cap \QQ$ which are arcs which are \bparallel\ on $P$ 
across halfdisks which do not contain $D_0$.

Now $\bd D_j = \bd D_j\p$ for a disk $D_j\p$ in $\bd \Yn0$ 
such that $D_j \cup D_j\p = \bd B_j$ for a 3-ball $B_j$ in $\Yn0$. 
Let $R\p$ be the union of $\Fn0$ and the $D_j\p$. 
Let $R=(\bd \Yn0)-int \, R\p$. Since each component of $\PP$ is 
non-compact we have that $\PP$ lies outside the union of the $B_j$ 
and thus $\PP \cap \QQ \cap \bd \Yn0$ lies in $R$. We claim that 
$R$ is \inc\ in $\Yn0$. For suppose $G$ is a compressing disk. 
Then $\bd G = \bd G\p$ for a disk $G\p$ in $\bd \Yn0$. Since 
no component of $\Fn0$ is a planar surface $G\p$ must lie in 
$(\bd \Yn0)-\Fn0$. But then $int \, G\p$ must contain some $D_j\p$, 
and so $D_j$ cannot meet $\Fn0$, a contradiction. 

Now let $\JJ$ be the union of all the simple closed curve components 
of $\PP \cap \QQ$. Then $\Yn0$, $\PP$, $\QQ$, and $\JJ$ satisfy the 
hypotheses of Proposition 2.1, and so $\JJ$ can be removed from the 
intersection by an isotopy which is fixed on those components of  
$\PP \cap \QQ$ which 
are not removed. Next let $\AAA$ be the union of all the components of 
$\PP \cap \QQ$ which are arcs that are \bparallel\ in $P$ across 
a halfdisk which does not contain $D_0$. Each such component is 
\bparallel\ in $\PP$ to an arc which does not meet $\Fn0$. Then 
$\Yn0$, $\PP$, $\QQ$, $R$, and $\AAA$ satisfy the hypotheses of 
Proposition 2.3, and so $\AAA$ can be removed from the intersection 
by an isotopy fixed on the components which are not removed.  
Note that since these isotopies are fixed on $R\p$ we may assume 
that they are fixed on the union of $\Fn0$ and the $B_j$. We thus 
get an ambient isotopy of $P$ in $M$ which puts it into \nzsp\ 
with respect to $C$. \qed \enddemo

Let $P$ be in \nzsp\ with respect to $C$; let $\{D_m\}$   
be the corresponding exhaustion for $P$. 
For each $m \geq 0$ each 
component $Z$ of $D_{m+1} - Int_P D_m$ is a partial disk 
and will be called a {\bf patch}. $Z$ lies in $X_{n+1}$ for some $n \geq n_0$ 
or in $C_{n_0} - Int \, C_0$. In the first case let 
$\bd_+ Z=\bd_1 Z \cap F_{n+1}$ and  $\bd_- Z=\bd_1 Z \cap F_n$; 
in the second case let $\bd_+ Z=\bd_1 Z$ and $\bd_- Z = \emptyset$. 
The set of all patches together with  $D_0$ forms the vertex set of a 
locally finite graph $\Gamma$ in  
which the edges are all the components of all the $Fr_P D_m$;   
two vertices are joined by an edge whenever the two disks meet along  
the corresponding arc. Since $P$ is simply connected $\Gamma$ 
is a tree. A subgraph of $\Gamma$ will often be identified with the 
submanifold of $P$ consisting of the union of its vertices. 
For $n \geq n_0$ let $\Gamma_n=\Gamma \cap C_n$. If $T$ is a component 
of $\Gamma_n$, then $T$ is a finite tree and $Fr_P T \subseteq F_n$. 
If $T$ does not contain $D_0$, then $T$ is called a {\bf falling tree} 
with {\bf frontier} in $F_n$. We also say that the falling tree $T$ 
{\bf descends} from $F_n$. A falling tree which is the star of some 
vertex $Z$ is called a {\bf falling star} with {\bf falling vertex} or 
{\bf center} $Z$. Every falling tree contains a falling star. 

$P$ is in {\bf \nzmp\ } with respect to $C$ if it is in 
\nzsp\ with respect to $C$ and for each patch $Z$ one 
has that $\bd_- Z$ is a single arc. This is equivalent to each 
$\Gamma_n$ being connected and hence to there being no falling trees. 

\proclaim{Lemma 4.2} Let $C$ be a nice exhaustion for $M$, and let $P$ be 
a proper partial plane in $M$. Then for some $n_0 > 0$ one has that 
$P$ is ambient isotopic to a partial plane which is in $n_0$-monotone 
position with respect to $C$. \endproclaim

\demo{Proof} Use Lemma 4.1 to put $P$ in \nzsp. By choosing $n_0$ 
sufficiently large we may assume that $D_0 \cap C_0 \neq \emptyset$. 
We first describe isotopies supported in 3-balls which reduce the number of 
edges of a falling star and which eliminate falling stars with two 
edges. We then examine the effect of these isotopies on other portions 
of $P$ and describe further isotopies which may be needed to keep $P$ 
in \nzsp. The concatenation of these isotopies gives an isotopy with 
compact support which eliminates a falling star but may create new ones. 
We then show how to organize these isotopies so as to eliminate all falling 
trees. 

Suppose $Z$ is the center of a falling star. Then $\bd_1 Z \subseteq 
F_{n+1}$ for some $n \geq n_0-1$, and $Z$ lies in $X$, where $X$ is  
$X_{n+1}$ if $n \geq n_0$ and $C_{n_0} - Int \, C_0$ if $n=n_0-1$. 
Let $X^+$ be the union of $X$ and a collar on $F_{n+1}$ in $X_{n+2}$. 
Since $X$ is \irr\ and \birr\ there is a disk $Z\p$ in $\bd X$ to which 
$Z$ is \bparallel\ across a 3-ball $B$ in $X$. Let $A$ be a component 
of $X \cap \bd M$ which meets $Z$. Then $A$ is an annulus which meets 
$F_{n+1}$ in a simple closed curve $K$. Since $Z \cap F_n = \emptyset$ 
each component $\be$ of $A \cap \bd Z$ is an arc which is \bparallel\ in 
$A$ to an arc $\be\p$ in $K$. Let $\be$ be an innermost such component, 
i.e. $\be \cup \be\p$ bounds a disk $G$ in $A$ such that $G \cap Z = \be$. 

Suppose $Z$ has order at least three. Let $B_1$ be a regular 
neighborhood of $G$ in $X^+$. There is an ambient isotopy of $P$ in 
$M$ supported in $B_1$ which moves $\be$ to $\be\p$ and then past 
$\be\p$ into $int(X_{n+2} \cap \bd M)$. Such an isotopy (regardless of  
the order of $Z$) is called a 
{\bf boundary  slide} of $\be$ past $\be\p$. The endpoints of $\be$ 
must lie on distinct arcs $\ga_1$, $\ga_2$ of $\bd_1 Z$; otherwise 
they would be joined by a single arc $\ga$ in $\bd_1 Z$ and so one 
would have $\bd Z = \be \cup \ga$, contradicting the fact that $P$ 
is in \nzsp. Now $\ga_1$ and $\ga_2$ are components of $\bd_- Z_1$ 
and $\bd_- Z_2$, respectively, for distinct patches $Z_1$ and $Z_2$ 
in $X_{n+2}$. The boundary slide replaces $Z_1$ and $Z_2$ by a new patch $W$ 
obtained by joining them by a band. $W$ has order one less than the 
sum of the orders of $Z_1$ and $Z_2$; since these orders are at least 
two $W$ has order at least three. The boundary slide replaces $Z$ by a partial 
disk $V$ with order one less than that of $Z$; hence $V$ has order 
at least two. 

Suppose $Z$ has order two. Then the endpoints of the arc 
$\ga_1 \cup \be \cup \ga_2$ are joined by an arc $\de$ in some component 
$A^*$ of $X \cap \bd M$, and the union of these two arcs is $\bd Z\p$. 
Moreover $\de$ is \bparallel\ in $A^*$ to an arc $\de\p$ in 
$A^* \cap F_{n+1}$ across a disk $G^*$. 

Suppose $G \cap G^* = \emptyset$. Then $\be\p \cup \ga_1 \cup \de\p 
\cup \ga_2$ bounds a disk $H$ in $F_{n+1}$, and 
$Z\p=G \cup H \cup G^*$. Let $B_2$ be a regular neighborhood in $X^+$ 
of $B$. There is an ambient isotopy of 
$P$ in $M$ supported in $B_2$ which carries $Z$ to $H$ and then past $H$ 
into $Int \, X_{n+2}$. This {\bf band push}  
replaces $Z$, $Z_1$, and $Z_2$ by a patch $Y$ in $X_{n+2}$ whose order 
is two less than the sum of the orders of $Z_1$ and $Z_2$ and hence 
is at least two. 

Suppose $G \cap G^* \neq \emptyset$. Then 
$A^*=A$, $G \subseteq Int_A G^*$, 
and $\be\p  \subseteq int \, \de\p$. Moreover each $\ga_i$ is 
\bparallel\ in 
$F_{n+1}$ across a disk $H_i$ to an arc $\de\p_i$ in $K$, $H_1 \cap H_2
= \emptyset$, $Z\p= H_1 \cup (G^* - Int_{G^*} G) \cup H_2$, and 
$\de\p = \de\p_1 \cup \be\p \cup \de\p_2$. Let $B_3$ be a regular  
neighborhood in $X^+$  of $B \cup G$. There is an ambient isotopy of 
$P$ in $M$ supported in $B_3$ which consists of a boundary slide of $\be$ 
past $\be\p$, followed by a boundary slide of $\de$ past $\de\p$,  
followed by a disk push which replaces $Z$, $Z_1$, and $Z_2$ by a 
patch $Y$ in $X_{n+2}$ of order at least two. This is a {\bf band 
unfolding}.

Now suppose $Z_0$ is some other vertex of $\Gamma \cap X$ which meets 
the ball $B_i$ ($i=1$, $2$, or $3$) supporting one of these isotopies. 
For a band push we have 
$Z_0 \subseteq B_2$, and so $Z_0$ is the center of a falling star. 
The isotopy replaces this star by a new vertex in $X_{n+2}$ of order at 
least two, and so $P$ is still in \nzsp. Now consider a boundary 
slide. Let $D=\be\p \times [0,1]$ be a disk in $F_{n+1}$ with 
$\be\p \times \{0\}=\be\p$, and $\be\p \times \{0,1\} \subseteq 
(\ga_1 \cup \ga_2)$. 
We may asssume  that $D \subseteq B_1$ and that $Z_0 \cap D$ consists 
of product arcs. This isotopy deletes each of the arcs of $Z_0 \cap (G \cup D)$ 
and joins its endpoints by an arc in $D$. It joins together the vertices in 
$X_{n+2}$ adjacent to $Z_0$ by bands corresponding to the new arcs in $D$ 
to get a 2-manifold $W_0$; it also replaces $Z_0$ by a disk $V_0$. 
Let $J=W_0 \cap V_0$. Then $J \subseteq \bd V_0$. If $J$ has at least two 
components, then each of these components is an arc, and $V_0$ and each 
component of $W_0$ is a partial disk of order at least two. If $J$ has 
exactly one component, then it is an arc or a simple closed curve. Suppose 
$J$ is an arc. Then $W_0$ is a partial disk of order at least three. 
If $V_0$ has order one, then $\bd V_0=\bd V_0\p$ for a disk $V_0\p$ in 
$\bd X$ such that $V_0\p$ is the union of a disk $V_1\p$ in $X \cap \bd M$ 
and a  disk $V_2\p$ in $F_{n+1}$ along an arc in $F_{n+1} \cap \bd M$. A 
halfdisk push of $V_0$ past $V_1\p$ across the ball $B_0$ bounded by 
$V_0 \cup V_0\p$ replaces $W_0$ and $V_0$ by a partial 
disk $W_1$ in $X_{n+2}$ of order at least two. If $Z_1$ is some other vertex 
in $X$ whose image under the boundary slide lies in $B_0$, then it must be 
the center of a falling star, and so the halfdisk push replaces the 
falling star by a vertex in $X_{n+2}$. Now suppose $J$ is a simple 
closed curve. Then $W_0$ is an annulus one of whose boundary components 
is $J=\bd V_0 \subseteq F_{n+1}$. So $V_0$ is \bparallel\ in $X$ to a disk 
$V_0\p$ in $F_{n+1}$ across a 3-ball $B_0$. A disk push of $V_0$ across 
$B_0$ past $V_0\p$ replaces $W_0$ and $V_0$ by a disk $W_1$ in $X_{n+2}$. 
Note that since $Z_0$ has order at least two $W_0$ must have been obtained 
by joining at least two distinct vertices in $X_{n+2}$. It follows that 
$W_1$ has order at least two. If $Z_1$ is some other vertex in $X$ whose 
image under the boundary slide lies in $B_0$, then again it must be the 
center of a falling star which is replaced by a vertex in $X_{n+2}$ by 
the disk push. Thus in all cases $P$ is put back in \nzsp. Finally 
consider a band unfolding. This move is the concatenation of a 
boundary slide and a halfdisk push. An analysis similar to that above 
provides isotopies which return $P$ to \nzsp\ and reduce the number 
of vertices in $\Gamma_{n+2}$.


Thus given any falling vertex $Z$ in $X$, there is an ambient isotopy
of $P$ in $M$ supported in $X^+$ which collapses the falling star with 
center $Z$ to a vertex in $X_{n+2}$. The only other possible effects 
of this isotopy on $\Gamma$ are to collapse other such falling stars 
with centers in $X$ in a similar fashion and to amalgamate vertices in 
$X_{n+2}$ which are adjacent to vertices in $X$, thereby reducing the orders 
of these vertices, but not reducing them to one. The number of vertices 
in each of $\Gamma \cap X$ and $\Gamma \cap X_{n+2}$ is reduced. 
$\Gamma$ is unchanged outside $X \cup X_{n+2}$. 

Define a sequence $n_0 < n_1 < n_2 < \cdots$ by choosing $n_{i+1}>n_i+1$ 
so that if $T$ is any falling tree descending from $F_{n_{i+1}}$, then 
$T \cap F_{n_i+1} = \emptyset$. Every falling tree in $C_{n_0} - Int \, C_0$ 
is a falling vertex. By the discussion above there is an isotopy 
supported in $C_{n_0+1}-Int \, C_0$ which eliminates them all. 
This isotopy may create new falling vertices descending from $F_{n_0+1}$, 
but it creates no new falling trees. It may reduce the order of $D_0$ but it 
does not eliminate it since $D_0 \cap C_0 \neq \emptyset$. 

Suppose $i>0$. Let $k$ be the smallest index for which there is a falling 
tree $T$ descending from $F_{n_i}$ such that $T \cap F_k \neq \emptyset$. 
Then $k>n_{i-1}+1$, and $T \cap X_k$ consists of falling vertices. 
There is an isotopy supported in $X_k \cup X_{k+1}$ which eliminates 
them and creates no new falling trees descending from $F_{n_i}$. 
Continuing in this fashion there is an isotopy supported in 
$Int \, C_{n_i+1}-C_{n_{i-1}+1}$ which eliminates all falling trees 
descending from $F_{n_i}$. Since these isotopies have disjoint compact 
supports they give a single ambient isotopy of $P$ in $M$ after which 
each $\Gamma_{n_i}$ is connected. 

Finally there is an isotopy supported in $(Int \, C_{n_i})-C_{n_{i-1}}$ 
which removes all falling trees descending from $F_n$ for 
$n_{i-1}<n<n_i$ and creates no new falling trees. Again this yields a 
single ambient isotopy of $P$ in $M$ after which $P$ is in 
monotone position. 
\qed  \enddemo

\proclaim{Theorem 4.3} Let $M$ be a connected, irreducible,  
non-compact 3-manifold which admits a nice exhaustion. 
Let $P$ be a proper partial plane in $M$. Then $\bd P$ lies in a 
single component $E$ of $\bd M$, and 
there is a collar $H$ on $E$ in $M$ such that $P \subseteq Int \, H$. 
\endproclaim

\demo{Proof} By Lemma 4.2 we may assume that $P$ is in 
$n_0$-monotone position 
with respect to a nice exhaustion $C$. 

Suppose $P$ meets at least two components of $\bd M$. Then for some 
$n$ there is a component $\al$ of $P \cap F_n$ which has one endpoint in 
a boundary plane $E$ and the other in another boundary plane $E\p$. 
Let $Z$ be the patch in $X_{n+1}$ containing 
$\al$. Then $\bd Z$ bounds a disk $Z\p$ in $\bd X_{n+1}$ by the 
$\bd$-irreducibility of $X_{n+1}$. But this is impossible because the simple 
closed curve $E \cap F_n$ meets $\bd Z$ transversely in a single point since 
$Z \cap F_n=\al$. Thus $\bd P$ lies in a single component $E$ of $\bd M$. 

We now claim that every component $\al$ of $P \cap F_n$ is parallel in $F_n$ 
to an arc $\be$ in $E \cap F_n$. For let $Z$ be the patch in $X_{n+1}$ 
containing $\al$. Then $\bd Z$ bounds a disk $Z\p$ in $\bd X_{n+1}$. 
The components of $F_n \cap E \cap Z\p$ are proper arcs in $Z\p$. If the 
endpoints of $\al$ lie on distinct such arcs $\be_1$ and $\be_2$, 
then the other endpoints of $\be_1$ and $\be_2$ must lie on components 
of $Z \cap F_n$ other than $\al$, contradicting the fact that $\bd_-Z=\al$. 
So the endpoints of $\al$ are joined by an arc $\be$. Then $\al \cup \be$ 
bounds a disk $W$ in $Z\p$. If $W$ does not lie in $F_n$, then $W$ must 
contain a component of $F_{n+1}$, contradicting the positive genus 
assumption on $C$. Thus $\al$ is parallel to $\be$ across the disk $W$ in 
$F_n$. 

For each component $\al$ of $P \cap F_n$ let $W(\al)$ be the disk in $F_n$ 
across which $\al$ is parallel to an arc $\be$ in $F_n \cap E$. Let 
$\Cal W_n$ be the set of maximal such disks, i.e. those $W(\al)$ which do 
not lie in $Int_{F_n} W(\al\p)$ for some $\al\p$. The $W(\al)$ are unique and 
the elements of $\Cal W_n$ 
are disjoint since otherwise some component of $F_n$ would be a disk. 
Let $G_n=E \cap C_n$ and $A_n = E \cap X_n$. 

Consider $\bd D_0$. It has not been assumed that $C_0$ is $\bd$-irreducible. 
Nevertheless, $\bd D_0$ does bound a disk $D^{\prime}_0$ in $\bd C_0$. 
For let $\Cal W_0=\{W(\al_{0,1}), \ldots , W(\al_{0,k_0})\}$. Then $G_0^+=
G_0 \cup W(\al_{0,1}) \cup \cdots \cup W(\al_{0,k_0})$ is a disk in $\bd C_0$ 
which contains $\bd D_0$, and the result follows. 

Now since $M$ is \irr\ $D_0$ is parallel across a 3-ball $B_0$ in $C_0$ to 
a disk $D^{\prime}_0$ in $\bd C_0$. There is a regular neighborhood $H_0$ of 
$G_0^+$ in $C_0$ such that $B_0 \subseteq Int \, H_0$. Then $G_0^*=
Fr_{C_0} H_0$ is a proper disk in $C_0$, and there is a product structure 
$G_0 \times [0,1]$ on $H_0$ such that $G_0 \times \{0\}=G_0$, 
$G_0 \times \{1\}=G_0^*$, and $L_0=(\bd G_0) \times [0,1] \subseteq F_0$. 

Let $\Cal W_1=\{W(\al_{1,1}), \ldots , W(\al_{1,k_1})\}$. 
Then $A_1^+=A_1 \cup L_0 \cup W(\al_{1,1}) \cup \cdots \cup W(\al_{1,k_1})$ 
is an annulus in $\bd X_1$ which contains $\bd Z$ for each patch $Z$ in 
$X_1$. Each such $Z$ is parallel across a 3-ball $B(Z)$ in $X_1$ to a 
disk $Z\p$ in $A_1^+$. There is a regular neighborhood $H_1$ of 
$A_1^+$ in $X_1$ such that each $B(Z)$ is contained in $Int \, H_1$ 
and $H_1 \cap F_0=H_0 \cap F_0$. Then $A_1^*=Fr_{X_1}H_1$ is a proper 
annulus in $X_1$, and there is a product structure $A_1 \times [0,1]$ 
on $H_1$ such that $A_1 \times \{0\}=A_1$, $A_1 \times \{1\}=A_1^*$, 
and $(\bd A_1) \times [0,1]=L_0 \cup L_1$, where $L_1$ is an annulus 
in $F_1$ and the product structures on $L_0$ induced by those on $H_0$ 
and $H_1$ agree. 

We now continue this process, constructing for each $n$ an 
$H_n=A_n \times [0,1]$ with $A_n \times \{0\}=A_n$, $A_n \times \{1\}=A_n^*$, 
a proper annulus in $X_n$, $(\bd A_n) \times [0,1]=L_{n-1} \cup L_n$, 
where $L_n$ is an annulus in $F_n$ and the product structures on $L_{n-1}$ 
induced by those on $H_{n-1}$ and $H_n$ agree such that $(P \cap X_n) 
\subseteq Int \, H_n$. We then let $H=\cup_{n \geq 0} H_n$. 
\qed \enddemo

\head 5. Strong Aplanarity: The General Case  \endhead

A {\bf partial plane system} is a surface $\Cal P$ each component of 
which is a partial plane. Thus an aplanar \3m\  $M$ is strongly aplanar 
if and only if for every proper partial plane system $\Cal P$ in $M$ 
there is a collar $H$ on $\bd M$ such that $Int \, H$ contains $\Cal P$. 
The goal of this section is to show that this is the case if $M$ 
is \irr\ and admits a nice exhaustion. 

\proclaim{Lemma 5.1} Suppose $M$ is a connected, irreducible, 
non-compact \3m\  which is not \homeo\  to $\halfspace$. Let $E$ be a 
component of $\bd M$ which is a plane and has the property that for 
every proper partial plane $P$ in $M$ with $\bd P$ in $E$ there is a 
collar $H$ on $E$ in $M$ such that $Int \, H$ contains $P$. Then for 
every proper partial plane system $\Cal P$ in $M$ with $\bd \Cal P$ 
in $E$ there is a collar $H\p$ on $E$ in $M$ such that $Int \, H\p$ 
contains $\Cal P$. \endproclaim

\demo{Proof} There is a collection $\{\be_{i,j}\}$ of disjoint arcs 
in $E$ whose union is end-proper in $E$ such that $\be_{i,j} \cap 
\Cal P = \bd \be_{i,j}$ and the union $K$ of $\Cal P$ with these arcs is 
simply connected. The notation is chosen so that $\be_{i,j}$ joins 
components $P_i$ and $P_j$ of $\Cal P$ for some choice of  $i<j$. 
 One way to see this is as follows. Let $\{D_n\}$ be 
an exhaustion for $E$ chosen so that each $D_n$ is a disk and 
$\bd \PP$ and $\cup \bd D_n$ are in general position. We may 
assume that $\bd D_0$ meets at least two components of $\PP$. 
There is then an arc $\be_{0,1}$ in $\bd D_0$ whose interior misses 
$\PP$ and whose endpoints lie in different components $P_0$ and 
$P_1$ of $\PP$. Let $K_1=P_0 \cup \be_{0,1} \cup P_1$. If $\bd D_0$ 
meets other components of $\PP$, then there is an arc $\be_{i,2}$ in 
$\bd D_0$ whose interior is disjoint from $K_1 \cup \PP$ such that 
one endpoint lies in a third component $P_2$ of $\PP$ and the other 
endpoint lies in $P_i$, where $i=0$ or $1$. If this endpoint of $\be_{i,2}$ 
is an endpoint of $\be_{0,1}$, then we isotop $\be_{i,2}$ slightly in $E$ 
to make the two arcs disjoint, keeping this endpoint of $\be_{i,2}$ in 
$P_i$. We then let $K_2=K_1 \cup \be_{i,2} \cup P_2$. We continue 
in this fashion until we have a simply connected 2-complex $K_{m_0}$ 
containing all those components of $\PP$ which meet $\bd D_0$. 
We then adjoin arcs in $\bd D_1$ and the components of $\PP$ which meet 
$\bd D_1$ but do not lie in $K_{m_0}$ to obtain a 2-complex $K_{m_1}$. 
Continuing in this manner we inductively construct the desired 2-complex $K$. 

Let $B_{i,j}$ be a regular neighborhood of $\be_{i,j}$ in the \3m\  obtained by 
splitting $M$ along $\Cal P$. Thus $B_{i,j}$ has the form 
$D_{i,j} \times [0,1]$, where $D_{i,j}$ is a halfdisk, $(\bd_0 D_{i,j}) 
\times [0,1]=B_{i,j} \cap E=B_{i,j} \cap \bd M$, $B_{i,j} \cap \Cal P 
=D_{i,j} \times \{0,1\}$, $B_{i,j} \cap P_i=D_{i,j} \times \{0\}$, and 
$B_{i,j} \cap P_j=D_{i,j} \times \{1\}$. Identify $D_{i,j}$ with 
$D_{i,j} \times \{1/2\}$. We now form the ``band sum''  $P$ of the 
components of $\PP$ along the $\be_{i,j}$ by deleting 
all the $Int_{\Cal P} (B_{i,j} \cap \Cal P)$ and then adjoining all 
the bands $(\bd_1 D_{i,j}) \times [0,1]$. Then $P$ is a proper partial 
plane in $M$ and so lies in a collar $H=E \times [0,1]$ on $E$ in $M$ 
with $E=E \times \{0\}$. Let $E^*=E \times \{1\}$. 

Let $M\p$ be the \3m\  obtained by splitting $M$ along $P$. Thus $\bd M\p$ 
contains two copies of $P$ whose identification recovers $M$. Let $M\pp$ 
be the component of $M\p$ containing $E^*$. Then $M\pp$ is a connected, 
\irr, non-compact \3m\  which is not \homeo\  to $\RRR$. Let $\Cal D$ be the 
union of those $D_{i,j}$ which lie in $M\pp$. Then $E^*$ and $\Cal D$ are 
proper \inc\  surfaces in $M\pp$ no component of which is a 2-sphere. 
If $E^*$ is trivial in $M\pp$, then it is trivial in $M$ and hence 
$M$ is \homeo\  to $\halfspace$, a contradiction. Therefore by Corollary 
2.2 there is an ambient isotopy of $E^*$ in $M\pp$, fixed on $\bd M\pp$, 
which takes $E^*$ to a plane disjoint from $\Cal D$. It follows that there  
is an ambient isotopy of $E^*$ in $M$, fixed on $P \cup \bd M$, which takes 
$E^*$ to a plane disjoint from all the $B_{i,j}$ and thus disjoint from 
$\Cal P$. The image $H\p$ of $H$ under this isotopy is a collar on $E$ 
in $M$ such that $Int \, H\p$ contains $\Cal P$. \qed \enddemo

\proclaim{Lemma 5.2} Suppose $M$ is a connected, \irr, non-compact  \3m\. 
Let $E$ be a component of $\bd M$ which is a plane. Let $H$ be a collar 
on $E$ in $M$. Suppose $\Cal P$ is a proper partial plane system in $M$ 
such that $\Cal P=\Cal P_0 \cup \Cal P_1$, where $\Cal P_0$ and $\Cal P_1$ 
are unions of components of $\Cal P$ such that $\Cal P_0 \subseteq H$ and 
$\Cal P_1 \cap E=\emptyset$. Then there is an ambient isotopy of $\Cal P_1$ 
in $M$, fixed on $\Cal P_0 \cup \bd M$, taking $\Cal P_1$ to a surface 
$\Cal P_1\p$ such that $\Cal P_1\p \cap H =\emptyset$. \endproclaim

\demo{Proof} We may assume that $\Cal P_1 \neq \emptyset$. Let $H=
E \times [0,1]$ with $E=E \times \{0\}$ and $E^*=E \times \{1\}$. 
Let $M\p$ be the \3m\  obtained by splitting $M$ along $\Cal P_0$. 
Let $M\pp$ be the component of $M\p$ containing $E^*$. Then  
$M\pp$ is a connected, \irr, non-compact \3m\  which is not \homeo\  to 
$\RRR$. $\Cal P_1$ and $E^*$ are proper \inc\  surfaces in $M\pp$ having 
no 2-sphere components. No component of $\Cal P_1$ is a plane. If 
$E^*$ is trivial in $M\pp$, then it is trivial in $M$ and so $M$ is 
\homeo\  to $\halfspace$. However, $\Cal P_1 \neq \emptyset$ and 
$\Cal P_1 \cap E = \emptyset$ implies that $\bd M$ has at least two 
components, a contradiction. Thus by Corollary 2.2 there is an ambient 
isotopy of $\Cal P_1$, fixed on $\bd M\pp$, which takes $\Cal P_1$ to 
a surface $\Cal P_1\p$ such that $\Cal P_1\p$ and $E^*$ are disjoint. 
The desired conclusion then follows. \qed \enddemo

\proclaim{Theorem 5.3} Let $M$ be a connected, \irr, 
non-compact \3m\  which admits a nice exhaustion. 
Then $M$ is strongly aplanar. \endproclaim

\demo{Proof} $M$ is aplanar by Theorem 3.4. 
$\bd M$ has components $E_1, \ldots, E_{\nu}$; each $E_i$ 
is a plane. Let $\Cal P$ be a proper partial plane system in $M$. 
By Theorem 4.3 $\Cal P=\Cal P_1 \cup \cdots \cup \Cal P_{\nu}$, 
where $\Cal P_i$ is the union of those components $P$ of $\Cal P$ such 
that $\bd P \subseteq E_i$. We induct on the number of nonempty $\Cal P_i$. 
If $\Cal P$ meets only one boundary component, say $E_1$, then we apply 
Theorem 4.3 and Lemma 5.1 to get a collar $H_1$ on $E_1$ in $M$ such 
that $\Cal P \subseteq Int \, H_1$ and then take arbitrary collars on 
the remaining $E_i$ in $M-Int \, H_1$. Suppose that the theorem is true 
for those proper partial plane systems meeting at most $k$ of the $\nu$ 
boundary components. Let $\Cal P$ be a proper partial plane system 
meeting $k+1$ of them, say $E_1, \ldots, E_{k+1}$. By Theorem 4.3 and 
Lemma 5.1 there is a collar $H_{k+1}=E_{k+1} \times [0,1]$ on $E_{k+1}$ 
in $M$ with $E_{k+1}=E_{k+1} \times \{0\}$ such that $\Cal P_{k+1} 
\subseteq Int \, H_{k+1}$. Let $M\p=M-Int \, H_{k+1}$. By Lemma 5.2 
there is an ambient isotopy of $\Cal P$ in $M$, fixed on $\Cal P_{k+1} 
\cup \bd M$, which takes $\Cal P_1 \cup \cdots \cup \Cal P_k$ to 
$\Cal P_1\p \cup \cdots \cup \Cal P_k\p$ lying in $M\p-(E_{k+1} \times 
\{1\})$. Since $M\p$ is \homeo\  to $M$ the inductive hypothesis gives 
a collar on $\bd M\p$ in $M\p$ whose interior in $M\p$ contains 
$\Cal P_1\p \cup \cdots \cup \Cal P_k\p$. It follows that there is a 
collar on $\bd M$ in $M$ whose interior in $M$ contains $\Cal P$. 
\qed \enddemo

\head 6. Attaching Boundary Planes  \endhead

In this section we show how to attach boundary planes to a connected, 
open \3m\  $U$ with finitely many ends to obtain non-compact 
3-manifolds $M$ with certain properties. In each case we construct $M$ 
by choosing a finite set of disjoint rays $[0,\infty)$ end-properly 
embedded in $U$ and then removing the interiors of disjoint regular 
neighborhoods of these rays. It is easily seen that $U$ is \homeo\  to 
$int \, M$ via a homeomorphism which takes the end of $U$ determined 
by a ray in $U$ to the end of $M$ determined by the corresponding 
boundary plane of $M$. The rays will be chosen using an exhaustion 
of $U$ so as to obtain a nice exhaustion for $M$. 

We will need some results and terminology from [14]. A compact, connected
\3m\  $X$ is {\bf excellent} if it is $\bold P^2$-\irr\ and \birr, it is not a 3-ball, 
it contains a 2-sided, proper, \inc\  surface, and every connected, 
proper, \inc\  surface of zero Euler characteristic in $X$ is \bparallel\. 
(So in particular $X$ is \ann\  and atoroidal.) Let $\la$ be a proper 
1-manifold in a \3m\  $Y$. The {\bf exterior} of $\la$ in 
$Y$ is the closure of the complement of a regular neighborhood of $\la$  
in $Y$. Suppose $Y$ is compact. If the exterior of $\la$ is excellent, 
then we say that $\la$ itself is {\bf excellent}.  It follows from Thurston's 
work [12] that $\la$ is excellent if and only if its exterior is hyperbolic. 
Theorem 1.1 of [14] states, among other things, that if $\kappa$ 
is any proper 1-manifold in a compact, connected \3m\  $Y$ such that 
$\kappa$ meets every 2-sphere in $\bd Y$ at least twice and every 
projective plane in $\bd Y$ at least once, then $\kappa$ is homotopic 
rel $\bd \kappa$ to an excellent 1-manifold $\la$. 

\proclaim{Theorem 6.1} Let $1 \leq \mu < \infty$, and for $1 \leq i \leq 
\mu$ let $1 \leq \nu_i < \infty$. Let $U$ be a connected, open  
\3m\  with $\mu$ ends. Then there is a non-compact \3m\ $M$ 
which has the following properties. 
\roster 
\item $U$ is \homeo\  to $int \, M$.
\item $\bd M$ is a disjoint union of planes $E^{i,j}$, 
$1\leq i\leq \mu$, $1\leq j \leq \nu_i$. 
\item For $1\leq i\leq \mu$ 
the image of the end $e^i$ of $U$ under the homeomorphism and 
inclusion induced maps $\varepsilon(U) \rightarrow \varepsilon(int \, M) 
\rightarrow \varepsilon(M)$ is 
the image of $\varepsilon(\cup_{j=1}^{\nu_i} E^{i,j})$ under the 
inclusion induced map $\varepsilon(\bd M) \rightarrow \varepsilon(M)$. 
\item $M$ is connected, eventually end-\irr, eventually $\bold P^2$-\irr, 
an\-an\-nu\-lar  at infinity, and totally acylindrical, and is not almost compact. 
\item If $U$ is \irr, then $M$ is strongly aplanar. 
\endroster
\endproclaim

\demo{Proof} Let $\{K_n\}$ be an exhaustion for $U$. Denote $M-K_n$ by 
$V_n$ and $K_{n+1}-Int \, K_n$ by $Y_{n+1}$. By passing to a subsequence of 
$\{K_n\}$ we may assume that each $V_n$ has $\mu$ components $V_n^i$,  
$1\leq i\leq \mu$, with $V^i_n \supseteq V^i_{n+1}$. Let $Y^i_{n+1}=
(Cl \, V^i_n) \cap Y_{n+1}$. By passing to a subsequence and attaching 
1-handles to $Fr \, K_n$ in $Cl \, V_n$ we may assume that each $Y^i_{n+1}$ 
is connected, each component of $Fr \, K_n$ has negative Euler 
characteristic, and each orientable component of $Fr \, K_n$ has positive 
genus. In particular $\bd Y_{n+1}$ contains no 2-spheres or projective 
planes. 

By Theorem 1.1 of [14] there exist disjoint proper arcs $\al^{i,j}_{n+1}$, 
$1\leq j\leq \nu_i$, in $Y^i_{n+1}$ such that $\al^{i,j}_{n+1}$ runs 
from $Fr \, K_n$ to $Fr \, K_{n+1}$, $\al^{i,j}_{n+1} \cap Fr \, K_{n+1}=
\al^{i,j}_{n+2} \cap Fr \, K_{n+1}$, and the exterior $X^i_{n+1}$ of 
$\cup_{j=1}^{\nu_i} \al^{i,j}_{n+1}$ in $Y^i_{n+1}$ is excellent. 
Then $\al^{i,j}=\cup_{n=1}^{\infty} \al^{i,j}_n$ is a ray which is proper 
in $V_0$ and end-proper in $U$. 

Let $X_{n+1}=\cup_{i=1}^{\mu} X^i_{n+1}$. We may assume that 
$X_{n+1}\cap Fr \, K_{n+1}=X_{n+2}\cap Fr \, K_{n+1}$. Let $M=K_0 \cup 
\cup_{n=1}^{\infty} X_n$. Note that $Cl_{Y_n}(Y_n-X_n)$ consists of 
disjoint 3-balls $N^{i,j}_n$ which are regular neighborhoods of 
$\al^{i,j}_n$, respectively. Let $N^{i,j}=\cup_{n=1}^{\infty} 
N^{i,j}_n$. Then $E^{i,j}=Fr_U N^{i,j}$ is a proper plane in $U$, and 
$\bd M$ is the disjoint union of these planes. 

Let $H^{i,j}$ be disjoint regular neighborhoods of $N^{i,j}$ in $U$. 
Then $H^{i,j}$ and $H^{i,j}-N^{i,j}$ are both \homeo\  to $\halfspace$, 
and thus there is a homeomorphism from $U$ to $int \, M$ with the 
required properties. 

Finally let $C_0=K_0$ and $C_{n+1}=C_n \cup X_{n+1}$. Then $\{C_n\}$ 
is a nice exhaustion for $M$, and so by Lemma 1.2, Lemma 1.3, Theorem 3.4, 
Theorem 3.5, and Theorem 5.3 the \3m\  $M$ has all the 
stated properties. \qed \enddemo

We now turn to the modifications of the basic construction of Theorem 6.1 
which yield the 3-manifolds with the additional properties mentioned 
in the introduction. We will need some preliminary results. 

\proclaim{Lemma 6.2} Let $R$ be a compact, connected \3m\. Let $S$ be a 
compact, proper, 2-sided surface in $R$. Let $R\p$ be the \3m\  obtained 
by splitting $R$ along $S$. Let $S_1$ and $S_2$ be the two copies of 
$S$ in $\bd R\p$ which are identified to obtain $R$. If each component 
of $R\p$ is excellent, $S_1 \cup S_2$ and $(\bd R\p)-int(S_1 \cup S_2)$ 
are \inc\  in $R\p$, and each component of $S$ has negative Euler 
characteristic, then $R$ is excellent. \endproclaim

\demo{Proof} This is Lemma 2.1 of [14]. \qed \enddemo

We first strengthen our construction so that it withstands the removal of 
some (but not all) boundary planes from each end; the basic technical 
device employed is called a ``poly-excellent tangle.'' 
Let $n\geq 1$. An {\bf $n$-tangle} is a proper 1-manifold $\la$ in a 
3-ball such that $\la$ has $n$ components and each of these components 
is an arc. For $1\leq k\leq n$ a {\bf $k$-subtangle} of an $n$-tangle 
$\la$ is a union of $k$ components of $\la$. An $n$-tangle $\la$ is 
{\bf poly-excellent} if for each $1\leq k\leq n$ each $k$-subtangle of 
$\la$ is excellent. 

\proclaim{Theorem 6.3} For all $n\geq 1$ poly-excellent $n$-tangles 
exist. \endproclaim

The proof of this theorem is given in the Appendix.

\proclaim{Lemma 6.4} Let $T$ be a solid torus or a solid Klein bottle, 
and let $1\leq \nu < \infty$. Then there exist disjoint proper arcs 
$\rho_1, \ldots, \rho_{\nu}$ in $T$ such that for every nonempty subset 
$\{j_1, \ldots, j_k\}$ of $\{1, \ldots, \nu\}$ 
the 1-manifold $\rho_{j_1} \cup \cdots \cup \rho_{j_k}$ is excellent 
in $T$. \endproclaim

\demo{Proof} Let $G$ be a meridional disk in $T$. Let $B$ be the 3-ball 
obtained by splitting $T$ along $G$, and let $G_1$ and $G_2$ be the disks 
in $\bd B$ which are identified to obtain $T$. 

By Theorem 6.3 $B$ contains a poly-excellent $3\nu$-tangle $\la$. Divide the 
components of $\la$ into three groups $\{\be_j\}$, $\{\ga_j\}$, and 
$\{\de_j\}$, where $1\leq j\leq \nu$. Isotop $\la$ so that $\be_j$ and 
$\ga_j$ run from $int \, G_1$ to $\bd B-(G_1 \cup G_2)$, $\de_j$ runs from 
$int \, G_2$ to itself, and under the identification of $G_1$ and $G_2$ one 
has that $(\be_j \cup \ga_j) \cap int \, G_1$ is identified with 
$\de_j \cap int \, G_2$. Let $\rho\p_j=\be_j \cup \ga_j \cup \de_j$, and 
let $\rho_j$ be the image of $\rho\p_j$ in $T$. 

Let $\{j_1, \ldots, j_k\}$ be a nonempty subset of $\{1, \ldots, \nu\}$. 
Let $R$ be the exterior of $\rho_{j_1} \cup \cdots \cup \rho_{j_k}$ in 
$T$, and let $R\p$ be the exterior of 
$\rho\p_{j_1} \cup \cdots \cup \rho\p_{j_k}$ in $B$. We may assume that 
these exteriors have been chosen so that $S=R \cap G$ is a disk with 
$2k$ holes and $R\p$ is the \3m\  obtained by splitting $R$ along $S$. 
Let $S_1$ and $S_2$ be the copies of $S$ in $\bd R\p$ which are identified 
to obtain $R$. Each component of $\bd R\p - int(S_1 \cup S_2)$ is an 
annulus. Since $\bd R\p$ is incompressible in $R\p$ and no component 
of $S_1 \cup S_2$ or of $\bd R\p-int(S_1 \cup S_2)$ is a disk each of 
these surfaces is incompressible in $R\p$. Since $R\p$ is excellent and 
$S$ has negative Euler characteristic Lemma 6.2 implies that $R$ is 
excellent. \qed \enddemo

\proclaim{Theorem 6.5} Let $\mu$, $\nu_i$, and $U$ be as in Theorem 6.1. 
Then there is a \3m\  $M$ having all the properties listed in that theorem 
such that if $1\leq \nu\p_i \leq \nu_i$, then every \3m\  $\widehat{M}$ 
obtained from $M$ by deleting for each $i$ any $\nu_i-\nu\p_i$ of the 
$E^{i,j}$ from $\bd M$ has all these properties. \endproclaim

\demo{Proof} Let $\{K_n\}$, $V_n$, $V^i_n$, $Y_{n+1}$, and $Y^i_{n+1}$ 
be as in the proof of Theorem 6.1. Let $\xi^i_{n+1}$ be a proper arc in 
$Y^i_{n+1}$ running from $Fr \, K_n$ to $Fr \, K_{n+1}$ chosen so that 
$\xi^i_{n+1} \cap Fr \, K_{n+1}=\xi^i_{n+2} \cap Fr \, K_{n+1}$. 
Let $Z^i_{n+1}$ be the exterior of $\xi^i_{n+1}$ in $Y^i_{n+1}$. 
We may assume that 
$Z^i_{n+1} \cap Fr \, K_{n+1}=Z^i_{n+2} \cap Fr \, K_{n+1}$. 
Then $\bd Z^i_{n+1}$ contains no 2-spheres or projective planes because 
each component of $Fr \, K_n \cup Fr \, K_{n+1}$ has negative Euler 
characteristic. By Theorem 1.1 of [14] there is an excellent proper arc 
$\eta^i_{n+1}$ in $Z^i_{n+1}$ with $\bd \eta^i_{n+1}$ in the open annulus 
$(\bd Z^i_{n+1}) \cap (int \, Y^i_{n+1})$. Let $L^i_{n+1}$ be the  
exterior of $\eta^i_{n+1}$ in $Z^i_{n+1}$. Then $T^i_{n+1}=
Y^i_{n+1}-Int_{Y^i_{n+1}} L^i_{n+1}$ is a solid torus or solid Klein 
bottle which meets $\bd Y^i_{n+1}$ in a disk $D^i_n$ in $Fr \, K_n$ 
and a disk $D^i_{n+1}$ in $Fr \, K_{n+1}$. 

By Lemma 6.4 there is a set of disjoint proper arcs $\al^{i,j}_{n+1}$,  
$1\leq j\leq \nu_i$, in $T^i_{n+1}$ with $\al^{i,j}_{n+1}$ running from 
$int \, D^i_n$ to $int \, D^i_{n+1}$ such that for every nonempty subset 
the union  of its elements is excellent in  $T^i_{n+1}$. We then proceed 
to construct $M$ as in the proof of Theorem 6.1, with the same notation 
as defined there. 

Now suppose that for each $i$ we have deleted a set of $\nu_i-\nu\p_i$ 
boundary planes $E^{i,j}$ from the $i^{th}$  end of $M$ so as to obtain 
$\widehat{M}$. This \3m\  is \homeo\  to the one obtained by ignoring the 
corresponding rays $\al^{i,j}$ in our construction. Hence this corresponds 
to ignoring the arcs $\al^{i,j}_{n+1}$. Let $\widehat{X}^i_{n+1}$ be 
the exterior of the remaining arcs in $Y^i_{n+1}$, and let $R^i_{n+1}$ be 
the exterior of these arcs in $T^i_{n+1}$. Then 
$\widehat{X}^i_{n+1}=L^i_{n+1} \cup R^i_{n+1}$ and the surface 
$S^i_{n+1}=L^i_{n+1} \cap R^i_{n+1}$ is a torus or Klein bottle from 
which the interiors of two disjoint disks have been removed. 
$\bd L^i_{n+1}-int \, S^i_{n+1}$ has two components obtained by removing 
$int \, D^i_n$ and $int \, D^i_{n+1}$ from components of $Fr \, K_n$ 
and $Fr \, K_{n+1}$ respectively. $\bd R^i_{n+1}-int \, S^i_{n+1}$ is 
a connected, orientable surface of genus $\nu\p_i-1$ with two boundary 
components. Thus $S^i_{n+1}$, $\bd R^i_{n+1}-int \, S^i_{n+1}$, and each 
component of  $\bd L^i_{n+1}-int \, S^i_{n+1}$ have negative Euler 
characteristic. It then follows from Lemma 6.2 that $\widehat{X}^i_{n+1}$ 
is excellent. We let $\widehat{X}_{n+1}=\cup_{i=1}^{\mu} 
\widehat{X}^i_{n+1}$, $\widehat{C}_0=C_0$, and $\widehat{C}_{n+1}=
\widehat{C}_n \cup \widehat{X}_{n+1}$. We then have a nice exhaustion 
$\{\widehat{C}_n\}$ of $\widehat{M}$. \qed \enddemo

We now further modify the construction so as to ensure the non-existence 
of homeomorphisms which carry one boundary plane to another or reverse 
orientation, as well as to produce an uncountable collection of pairwise 
non-\homeo\  3-manifolds having the same given interior and distribution of 
boundary planes among the ends. 

A {\bf classical knot space} $Q$ is a space \homeo\  to the exterior of 
a non-trivial knot in $S^3$. If $Q$ is embedded in a \3m\  $X$ and $\bd Q$ 
is \inc\  in $X$, then we say that $Q$ is {\bf incompressibly embedded} 
in $X$. The idea of the proof is to associate to each boundary plane $E$ an 
infinite sequence of disjoint classical knot spaces whose union lies in $\inte M$ 
and is 
end-properly embedded in $M$. These knot spaces are not incompressibly 
embedded in $M$. However, for any ``sufficiently large'' compact subset $K$ 
of $M$ which meets $E$ all but finitely many of them will be incompressibly 
embedded in $M-K$. On the other hand at most finitely many will be 
incompressibly 
embedded in the complement of a compact subset which does not meet $E$. 
Furthermore the knot spaces will be chosen so that they characterize the 
plane with which they are associated. A useful analogy is that of a string of 
lights associated to each plane, with different planes corresponding to 
different colors. The removal of certain compact subsets meeting a collection 
of planes then turns on all 
but finitely many lights in the strings associated to these planes.  

We first need a couple of preliminary technical lemmas. 

\proclaim{Lemma 6.6} Let $T$ be a solid torus or solid Klein bottle,  
and let $1\leq \nu < \infty$. Then there exist disjoint solid tori 
$T_1, \ldots, T_{\nu}$ 
in $int \, T$, disjoint compressing disks $D_1, \ldots, D_{\nu}$ for 
$\bd T_1, \ldots \bd T_{\nu}$ in $T-int(T_1 \cup \cdots \cup T_{\nu})$, 
and disjoint proper arcs $\rho_1, \ldots, \rho_{\nu}$ in $T-(T_1 \cup \cdots 
\cup T_{\nu})$ such that $D_i \cap \rho_j = \emptyset$ for $i\neq j$ 
and for every nonempty subset $\{j_1, \ldots, j_k\}$ 
of $\{1, \ldots, \nu\}$ the 1-manifold 
$\rho_{j_1} \cup \cdots \cup \rho_{j_k}$ is excellent in 
$T-int(T_{j_1} \cup \cdots \cup T_{j_k})$. \endproclaim

\demo{Proof} Let $G$, $B$, $G_1$, and $G_2$ be as in the proof of 
Lemma 6.4. 
Let $W_1, \ldots, W_{\nu}$ be disjoint disks in $int \, G$. 
Let $W_{i,j}$ be the copy of $W_j$ in $G_i$ whose image under the 
identification of $G_1$ with $G_2$ is $W_j$. 
Let $T_1, \ldots, T_{\nu}$ be disjoint regular neighborhoods of 
$\bd W_1, \ldots, \bd W_{\nu}$ in $T$, chosen so that $A_j=T_j \cap G$ is 
a regular neighborhood of $\bd W_j$ in $G$. Then $T_j$ is split 
by $A_j$ into solid tori $T_{1,j}$ and $T_{2,j}$ which are regular 
neighborhoods of $\bd W_{1,j}$ and $W_{2,j}$, respectively, in $B$. 
Let $A_{i,j}=Fr_B T_{i,j}$. Let $T^*=T-int(T_1 \cup \cdots \cup T_{\nu})$ and 
$B^*=B-Int_B (T_{1,1} \cup \cdots \cup T_{1,\nu} \cup T_{2,1} \cup 
\cdots \cup T_{2,\nu})$. Then $B^*$ is a 3-ball which meets $G_i$ in 
$\nu$ disks $D_{i,j} \subseteq int \, W_{i,j}$ and a disk with $\nu$ 
holes $H_i$. Let $D_j$ be the image of $D_{i,j}$ under the identification. 
Then the $D_j$ are disjoint compressing disks for $\bd T_j$ in $T^*$. 

By Theorem 6.3  $B^*$ contains a poly-excellent $4\nu$-tangle $\la$. 
Divide the components of $\la$ into four groups $\{\be_j\}$,  $\{\ga_j\}$, 
$\{\de_j\}$, and $\{\omega_j\}$, $1\leq j\leq \nu$. Isotop $\la$ so 
that $\be_j$ runs from $\bd B-(G_1 \cup G_2)$ to $int \, D_{1,j}$, 
$\ga_j$ runs from $int \, D_{2,j}$ to itself, $\de_j$ runs from 
$int \, D_{1,j}$ to $int \, H_1$, and $\omega_j$ runs from $int \, H_2$ to 
$\bd B-(G_1 \cup G_2)$. Do this so that under the identification we 
have $(\be_j \cup \de_j) \cap D_{1,j}$ identified with $\ga_j \cap D_{2,j}$ 
and we have $\de_j \cap H_1$ identified with $\omega_j \cap H_2$. 
Let $\rho\p_j=\be_j \cup \ga_j \cup \de_j \cup \omega_j$. Let $\rho_j$ 
be the image of $\rho\p_j$ under the identification. 

Now let $\{j_1, \dots, j_k\}$ be a nonempty subset of $\{1, \ldots, \nu\}$. 
Let $R_0$ be the exterior of $\rho_{j_1} \cup \cdots \cup \rho_{j_k}$ 
in $T-int(T_{j_1} \cup \cdots \cup T_{j_k})$, $R_0\p$ the exterior of 
$\rho\p_{j_1} \cup \cdots \cup \rho\p_{j_k}$ in $B$, and $R^*_0$ the 
exterior of $\rho\p_{j_1} \cup \cdots \cup \rho\p_{j_k}$ in $B^*$. 
We may assume these exteriors are chosen so that $R_0\p$ is the union of 
$R^*_0$ with all those $T_{i,j}$ for which $j \in \{j_1, \ldots, 
j_k\}$. Since $A_j$ is parallel to $A_{i,j}$ across $T_{i,j}$ we have 
that $R^*_0$ is \homeo\ to $R\p_0$ and is therefore excellent. We may 
further assume that $R^*_0$ is obtained by splitting $R_0$ along the 
surface $S$ consisting of the $k$ disks with two holes $R_0 \cap D_{j_r}$ 
and the disk with $2k$ holes $R_0 \cap (G-(D_{j_1} \cup \cdots \cup 
D_{j_k}))$. Let $S_1$ and $S_2$ be the two copies of $S$ in $\bd R^*_0$. 
Then the components of $\bd R^*_0-int(S_1 \cup S_2)$ are $4k$ annuli 
and a disk with $2k+1$ holes. Since $\bd R^*_0$ is \inc\  in $R^*_0$ 
it follows that $S_1 \cup S_2$ and $\bd R^*_0-int(S_1 \cup S_2)$ are 
\inc\  in $R^*_0$. Thus by Lemma 6.2 we have that $R_0$ is excellent. 
\qed \enddemo

\proclaim{Lemma 6.7} Let $T$ be a solid torus or solid Klein bottle. 
Let $J_1, \ldots, J_{\nu}$ be excellent knots in $S^3$. Then there are 
disjoint classical knot spaces $Q_1, \ldots, Q_{\nu}$ in $int \, T$ 
and disjoint proper arcs $\rho_1, \ldots, \rho_{\nu}$ in 
$T-int(Q_1 \cup \cdots \cup Q_{\nu})$ such that  
$Q_j$ is \homeo\  to the exterior of $J_j$ in $S^3$, 
there are disjoint 
3-balls $B_j$ in $int \, T$ such that $Q_j \subseteq B_j$ and 
$B_i \cap \rho_j=\emptyset$ for $i \neq j$. 
and for every nonempty subset 
$\{j_1, \ldots, j_k\}$ of $\{1, \ldots, \nu\}$ 
\roster
\item  the exterior $R$ of the 
1-manifold $\rho_{j_1} \cup \cdots \cup \rho_{j_k}$ in $T$ is 
$\bold P^2$-\irr, $\bd$-irreducible, and \ann, 
\item each $Q_{j_r}$ is 
incompressibly embedded in $R$, and 
\item given any classical knot space $Q$ 
incompressibly embedded in $int \, R$ there is an ambient isotopy of 
$Q$ in $R$, fixed on $\bd R$, which takes $Q$ to some $Q_{j_r}$. 
\endroster 
\endproclaim

\demo{Proof} Let $T_j$, $D_j$, and $\rho_j$ be as in Lemma 6.6. Let 
$T^*=T-int(T_1 \cup \cdots \cup T_{\nu})$. Let $Q_j$ be the exterior of 
$J_j$ in $S^3$. Form $T_0$ by gluing $Q_1 \cup \cdots \cup Q_{\nu}$ to 
$T^*$ by identifying $\bd Q_j$ with $\bd T_j$ so that a meridian of 
$J_j$ is identified with $\bd D_j$. The union $B_j$ of $Q_j$ with a 
regular neighborhood of $D_j$ in $T^*$ is then a 3-ball, and $T_0$ is 
again a solid torus or solid Klein bottle. 

Let $R_0$ be the exterior of $\rho_{j_1} \cup \cdots \cup \rho_{j_k}$ 
in $T-int(T_{j_1} \cup \cdots T_{j_k})$. Then $R=R_0 \cup Q_{j_1} \cup 
\cdots \cup Q_{j_k}$ is the exterior of  $\rho_{j_1} \cup \cdots \cup 
\rho_{j_k}$ in $T_0$. Since $R_0$ and the $Q_{j_r}$ are $\bold P^2$-\irr, \birr, 
\ann, and atoroidal and $\bd R$ is not a torus one can apply standard 
general position 
and isotopy arguments to show that each $\bd Q_{j_r}$ is \inc\  in $R$ and 
that every  \inc\  torus in $R$ is isotopic to some $\bd Q_{j_r}$. 
The result follows after changing the name of $T_0$ to $T$. 
\qed \enddemo

\proclaim{Theorem 6.8} Let $\mu$, $\nu_i$, and $U$ be as in Theorem 6.5. 
\roster 
\item There is a \3m\  $M$ having all the properties listed in Theorem 6.5 
such that $M$ admits no self homeomorphisms which take one boundary 
plane to another or reverse orientation. 
\item Each $\widehat{M}$ as in 
Theorem 6.5 also has all these properties (including (1)), and distinct $\widehat{M}$ 
are not \homeo. 
\item There are uncountably many pairwise non-homeomorphic 
such $M$, and if $M$ and $N$ are two of these manifolds which are 
not homeomorphic, then for every pair of associated manifolds $\widehat{M}$ 
and $\widehat{N}$ we have that $\widehat{M}$ and $\widehat{N}$ are not 
\homeo. \endroster \endproclaim 

\demo{Proof} Let $\Cal J$ be a countably infinite set of excellent 
knots in $S^3$ whose exteriors are pairwise non-homeomorphic and admit 
no orientation reversing self homeomorphisms. An example of such a set 
is the collection of all non-trivial twist knots other than the trefoil 
and figure eight knots [16]. Let $\Cal S$ be the set of all triples 
$(i,j,n)$ where $1\leq i\leq \mu$, $1\leq j\leq \nu_i$, and $n\geq 1$. 
Index $\Cal J$ by choosing a bijection with $\Cal S \times \{0,1\}$. 
Denote the indexed knot by $J(i,j,n,p)$, where $p \in \{0,1\}$, and 
its exterior by $Q(i,j,n,p)$. 

Let $\varphi:\Cal S \rightarrow \{0,1\}$. We will associate a \3m\  
$M$ to $\varphi$ as follows. Let $K_n$, $L^i_{n+1}$, $T^i_{n+1}$, 
and $D^i_n$ be as in the proof of Theorem 6.5. We apply Lemma 6.7 
to $T^i_{n+1}$ together with the knots $J(i,j,n+1,\varphi(i,j,n+1))$, 
$1\leq j\leq \nu_i$ to get disjoint proper arcs $\al^{i,j}_{n+1}$ running 
from $D^i_n$ to $D^i_{n+1}$ having the properties stated for the 
$\rho_j$ in the lemma. We then carry out the rest of the construction 
in the proof of Theorem 6.5 to get $M$. 

Suppose $\widehat{M}$ is obtained by deleting boundary planes 
as in Theorem 6.5. Then $\widehat{M}$ has an exhaustion 
$\{\widehat{C}_n\}$ with $\widehat{V}^i_n$ the $i^{th}$ component of 
$\widehat{M}-\widehat{C}_n$, $\widehat{V}^i_n=\cup_{q=n}^{\infty} 
\widehat{X}^i_{q+1}$, and $\widehat{X}^i_{q+1}=L^i_{q+1} \cup R^i_{q+1}$, 
where $R^i_{q+1}$ is the exterior of $\al^{i,j_1}_{q+1} \cup \cdots 
\cup \al^{i,j_k}_{q+1}$ in $T^i_{q+1}$. 

Standard general position and isotopy arguments show that 
$\widehat{X}^i_{q+1}=L^i_{q+1} \cup R^i_{q+1}$ is $\bold P^2$-\irr, \birr, and 
\ann. It follows that $\{\widehat{C}_n\}$ is a nice exhaustion for  
$\widehat{M}$. Arguments of this type also show that each $Q(i,j_r, q+1, 
\varphi(i,j_r, q+1))$ in our construction is incompressibly embedded in 
$\widehat{V}^i_n$ whenever $q\geq n$ and that every classical knot space 
$Q$ which is incompressibly embedded in $\widehat{V}^i_n$ can be 
ambiently isotoped into $\widehat{X}^i_{q+1}$ for some $q\geq n$, then 
into $R^i_{q+1}$, and hence by Lemma 6.7 to some 
$Q(i,j_r, q+1, \varphi(i,j_r,q+1))$ in our construction. 

Now suppose we have another function $\psi:\Cal S \rightarrow \{0,1\}$. 
Denote the two resulting $M$ by $M[\varphi]$ and $M[\psi]$, and 
distinguish the various submanifolds arising in their construction by 
similar notation. Let $\widehat{M}[\psi]$ be obtained by deleting all 
but one boundary plane from each end of $M[\psi]$. Suppose 
$g:\widehat{M}[\psi] \rightarrow N$ is a homeomorphism, where $N$ is 
obtained by deleting some boundary planes from $M[\varphi]$. Since $g$ 
induces a bijection $\varepsilon(\widehat{M}[\psi]) \rightarrow 
\varepsilon(N)$ and $\widehat{M}[\psi]$ has exactly one boundary 
plane per end, so does $N$. Hence it must be some $\widehat{M}[\varphi]$ 
obtained by deleting all but one boundary plane from each end of 
$M[\varphi]$. 

Fix $i$, and let $E^{i,j}[\psi]$ be the single boundary plane of the 
$i^{th}$ end of $\widehat{M}[\psi]$. Then $g(E^{i,j}[\psi])=
E^{s,t}[\varphi]$ for some $1\leq s\leq \mu$ and $1\leq t\leq \nu_s$. 
Choose $n>0$ such that $g(\widehat{C}_0[\psi]) \subseteq 
Int \, \widehat{C}_n [\varphi]$. Then choose $m>0$ such that 
$\widehat{C}_n [\varphi] \subseteq Int \, g(\widehat{C}_m [\psi])$. 
Let $q\geq m$, and let $Q$ be the copy of $Q(i,j,q+1,\psi(i,j,q+1))$ 
embedded in $\widehat{M}[\psi]$ by our construction. 
Then $Q$ lies in $\widehat{V}^i_m[\psi]$, and $\bd Q$ is \inc\  in 
$\widehat{V}^i_0 [\psi]$. Therefore $g(Q)$ lies in $g(\widehat{V}^i_m 
[\psi])$ and thus in the larger set $\widehat{V}^s_n [\varphi]$, and 
$g(\bd Q)$ in \inc\  in $g(\widehat{V}^i_0 [\psi])$ and thus in the 
smaller set $\widehat{V}^s_n [\varphi]$. Hence $g(Q)$ is isotopic to 
the copy of some $Q(s,t,r+1, \varphi(s,t,r+1))$ embedded in 
$\widehat{M}[\varphi]$ by our construction. Therefore $s=i$, $t=j$, 
$r=q$, and $\varphi(i,j,q+1)=\varphi(i,j,q+1)$. 

Now suppose $h:M[\psi] \rightarrow M[\varphi]$ is a homeomorphism. 
By restricting $h$ to $\widehat{M}[\psi]$ as above we see that $h$ 
must take the $i^{th}$ end of $M[\psi]$ to the $i^{th}$ end of 
$M[\varphi]$ and the $j^{th}$ boundary plane of the $i^{th}$ end of 
$M[\psi]$ to the $j^{th}$ boundary plane of the $i^{th}$ end of 
$M[\varphi]$. Moreover, there exists $m>0$ such that $\psi(i,j,q+1)=
\varphi(i,j,q+1)$ for all $q \geq m$. Thus for fixed $i$ and $j$ we get two 
infinite sequences of zeros and ones which agree after a finite number 
of terms. This property defines an equivalence relation on the set of 
all such sequences, which is uncountable, such that each equivalence 
class is countable. Therefore the set of all equivalence classes is 
uncountable, and so the set of homeomorphism classes of the $M[\varphi]$ 
is uncountable. 

Taking $\psi=\varphi$ we get that $h$ must take each boundary plane 
to itself and some $Q$ to itself. Since these classical knot spaces 
admit no orientation reversing homeomorphisms neither does $M[\varphi]$. 

Clearly these considerations apply to all the $\widehat{M}[\psi]$ and 
$\widehat{M}[\varphi]$ obtained by deleting boundary planes as in 
Theorem 6.5, and so we are done. \qed \enddemo

\head Appendix: Poly-excellent Tangles  \endhead

Recall that an $n$ component tangle in a 3-ball is poly-excellent if for 
every $k$ with $1 \leq k \leq n$ each of its $k$ component subtangles is 
excellent, 
i.e\. has  hyperbolic exterior. In this appendix we prove Theorem 6.3, 
which asserts the existence of poly-excellent $n$-tangles for all $n \geq 1$. 
The case $n=1$ is trivial since we can choose any proper arc in a 
3-ball having exterior homeomorphic to the exterior of an excellent knot 
in $S^3$, for example the figure eight knot. So we may assume $n \geq 2$. 

We shall make use of an excellent $n$-tangle $\la$ in a 3-ball $B$, 
defined for $n \geq 2$, called the true lover's $n$-tangle. This 
$n$-tangle is defined and its basic properties are proven in section 4 
(pages 275-281) of [13]. See figures 1 and 2 of that paper. Each component 
$\la_j$ of $\la$ is a trefoil knotted arc. Two distinct components 
$\la_j$ and $\la_i$ are linked in $B$ if and only if $|j-i|=1$. In fact, 
$B=B_1 \cup \cdots \cup B_{2n-1}$, where each $B_p$ is a 3-ball, 
$B_p \cap B_{p+1}$ is a disk, $B_p \cap B_q = \emptyset$ for $|p-q| > 1$, 
$\la_1 \subseteq B_1 \cup B_2$, 
$\la_j \subseteq B_{2j-2} \cup B_{2j-1} \cup B_{2j}$ for $1 < j < n$, 
and $\la_n \subseteq B_{2n-2} \cup B_{2n-1}$. 
Moreover $\bd \la_j \subseteq \bd B_{2j-1}$ for all $j$. See figure 3 of [13]. 

Now $\la$ is excellent (Proposition 4.1 of [13]) but is not poly-excellent. 
However, it has the property that for $2 \leq k \leq n$ each $k$-tangle 
consisting of $k$ consecutive components of $\la$ is excellent. This 
property is not stated explicitly in [13], but it follows directly from 
the proof of Proposition 4.1 of that paper. 

The basic idea behind the proof of Theorem 6.3 is to stack up several 
copies of $\la$, joining the bottom endpoints of one copy to the top 
endpoints of the copy beneath it, to obtain a new $n$-tangle $\theta$. 
The endpoints are to be joined using braids, so that each component of 
$\theta$ consists of segments which are components of the copies of $\la$ 
and may have different indices. This is to be done so that given any subset 
of $\{1, \ldots, n\}$ any two components of $\theta$ with indices in the 
subset will have segments with consecutive indices in some copy of $\la$. 
The exterior of the resulting subtangle $\widehat{\theta}$ of $\theta$ is 
then to be analyzed using Lemma 6.2. 

Unfortunately, this basic idea does not work. This can be seen by 
considering the case where $\widehat{\theta}$ has exactly one component. 
It then meets each disk between adjacent 3-balls in a single point, from 
which it follows that its exterior in not anannular. It also meets each 
of these 3-balls in a trefoil knotted arc, which is not excellent. 

Fortunately, there is a modification of the basic idea which does work. The 
$n$-tangle $\theta$ will be constructed so that, among other things, each 
of its components doubles back twice at each level so that it meets each 
intermediate disk three times and meets each 
3-ball in an excellent tangle. This requires us to use copies of the true 
lover's tangle having more than $n$ components, which accounts for much of 
the complication in the following argument. 

\demo{Proof of Theorem 6.3} As noted above we may assume $n \geq 2$. 
We take as the 3-ball containing $\th$ the set 
$$B=\{(x,y,z): 0 \leq x \leq 9n+1, -1 \leq y \leq 1, 0 \leq z \leq n^2-n+1\}.$$ 
We regard $x$ and $y$ as increasing from left to right and from 
back to front, respectively, and $z$ as increasing in the 
downward direction. For $0 \leq p \leq n^2-n+1$, $1 \leq q \leq 3n$, 
and $1 \leq j \leq n$, let $H_p=[0,9n+1] \times [-1,1] \times \{p\}$, 
$x_{p,q}=(3q-1,0,p)$, $a_{p,j}=x_{p,3j-2}$, $b_{p,j}=x_{p,3j-1}$, and 
$c_{p,j}=x_{p,3j}$. 

Let $m=(n^2-n)/2$. We now define certain subsets of $B$. 

First suppose $0 \leq i \leq m$. 
$$ B_i=[0,9n+1] \times [-1,1] \times [2i,2i+1] $$
$$ B_{i,j}=[9j-9,9j+1] \times [-1,1] \times [2i,2i+1], 1 \leq j \leq n $$
$$ N_{i,j}=[9j,9j+1] \times [-1,1] \times [2i,2i+1], 0 \leq j \leq n $$

Next suppose $1 \leq i \leq m$. 
$$C_i=[0,9n+1] \times [-1,1] \times [2i-1,2i] $$ 
$$C_{i,j}=[9j-9,9j+1] \times [-1,1] \times [2i-1,2i], 1 \leq j \leq n $$ 
$$K_{i,j}=[9j,9j+1] \times [-1,1] \times [2i-1,2i], 0 \leq j \leq n $$ 

Thus $B$ is a stack of rectangular solids, starting with $B_0$ on the top and then 
alternating in the pattern $C_i$, $B_i$, $C_{i+1}$, $B_{i+1}$ until 
concluding with $B_m$ on the bottom. We have $C_i \cap B_i=H_{2i}$ and 
$B_i \cap C_{i+1}=H_{2i+1}$. Each $B_i$ consists of $n$ rectangular solids  
$B_{i,j}$ which are disjoint except for the pairs $B_{i,j}$ and $B_{i,j+1}$, 
whose overlap is the rectangular solid $N_{i,j}$. 
A similar pattern holds for the $C_{i,j}$ and $K_{i,j}$. We will refer to 
the $B_{i,j}$ and $C_{i,j}$ as {\bf blocks}. Let $N_i=N_{i,0} \cup \cdots 
\cup N_{i,n}$ and $K_i=K_{i,0} \cup \cdots \cup K_{i,n}$. Let 
$B\p_{i,j}$ be the closure in $B$ of $B_{i,j}-N_i$. We define $C\p_{i,j}$ 
in a similar way. The $B\p_{i,j}$ and $C\p_{i,j}$ are called {\bf bricks}. 

Let $\La_0$ be a copy of the true lover's $2n$-tangle in $B_0$ with 
components $\la_{0,1}$, $\ldots$, $\la_{0,2n}$. For $1 \leq j \leq n$ 
let $\al_{0,j}=\la_{0,2j-1}$, $\ga_{0,j}=\la_{0,2j}$, and $\La_{0,j}=
\al_{0,j} \cup \ga_{0,j}$. Isotop $\La_0$ so that $\La_{0,j}$ is a 
$2$-tangle in $B_{0,j}$, $\al_{0,j}$ runs from $a_{0,j}$ to $a_{1,j}$, 
and $\ga_{0,j}$ runs from $b_{1,j}$ to $c_{1,j}$. The existence of this 
isotopy and the similar isotopies required below follows from the 
description of $\la$ given earlier in this appendix. 

For $1 \leq i \leq m-1$ let $\La_i$ be a copy of the true lover's 
$3n$-tangle in $B_i$ with components $\la_{i,1}$, $\ldots$, $\la_{i,3n}$. 
For $1 \leq j \leq n$ let $\de_{i,j}=\la_{i,3j-2}$, $\al_{i,j}=\la_{i,3j-1}$, 
$\ga_{i,j}=\la_{i,3j}$, and $\La_{i,j}=\de_{i,j} \cup \al_{i,j} \cup 
\ga_{i,j}$. Isotop $\La_i$ so that $\La_{i,j}$ is a $3$-tangle in $B_{i,j}$, 
$\de_{i,j}$ runs from $a_{2i,j}$ to $b_{2i,j}$, $\al_{i,j}$ runs from 
$c_{2i,j}$ to $a_{2i+1,j}$, and $\ga_{i,j}$ runs from $b_{2i+1,j}$ to 
$c_{2i+1,j}$. (For $n=2$ this piece of the construction does not occur.) 

Let $\La_m$ be a copy of the true lover's $2n$-tangle in $B_m$ with 
components $\la_{m,1}$, $\ldots$, $\la_{m,2n}$. For $1 \leq j \leq n$ 
let $\de_{m,j}=\la_{m,2j-1}$, $\al_{m,j}=\la_{m,2j}$, and $\La_{m,j}=
\de_{m,j} \cup \al_{m,j}$. Isotop $\La_m$ so that $\La_{m,j}$ is a 
$2$-tangle in $B_{m,j}$, $\de_{m,j}$ runs from $a_{2m,j}$ to $b_{2m,j}$, 
and $\al_{m,j}$ runs from $c_{2m,j}$ to $c_{2m+1,j}$. 

Let $\Cal B_{3n}$ be the Artin braid group on $3n$ strings; let 
$\si_1, \ldots, \si_{3n-1}$ be the standard generators for $\Cal B_{3n}$. 
(See [1].) For $1 \leq i \leq m$, given an element $\beta_i$ of 
$\Cal B_{3n}$, we interpret it as a geometric braid in $C_i$, i.e\. 
it consists of $3n$ disjoint proper arcs in $C_i$ such that the 
$q^{th}$ arc runs from $x_{2i-1,q}$ to some $x_{2i,r}$ and meets each 
horizontal plane in a single point. We follow the convention that 
as one reads a word in the generators of $\Cal B_{3n}$ from left 
to right the geometric braid goes downward. 

Thus we can associate to each sequence $\be_1, \ldots, \be_m$ of 
elements of $\Cal B_{3n}$ a proper 1-manifold $\th$ in $B$ by taking the 
union of the $\be_i$, $1 \leq i \leq m$, and the $\La_i$, $0 \leq i \leq m$. 
For $1 \leq j \leq n-1$, let 
$$ \Si_j = \si_{3j} \si_{3j-1} \si_{3j+1} \si_{3j-2} \si_{3j} \si_{3j+2} 
\si_{3j-1} \si_{3j+1} \si_{3j}.$$ 
If one partitions the $3n$ strings into $n$ consecutive groups of $3$ 
consecutive strings, then $\Si_j$ is obtained by crossing the $j^{th}$ 
group in front of the $(j+1)^{st}$ group and is thus the analogue of 
the $j^{th}$ standard generator of $\Cal B_n$. Note that if $\be_i=\Si_j$, 
then we may assume that the strings numbered $3j-2$ through $3j+3$ lie 
in $C_{i,j} \cup C_{i,j+1}$ and that all other strings are vertical. 

We now let the sequence $\be_1, \ldots, \be_m$ be 
$$\Si_1, \ldots, \Si_{n-1}, \Si_1, \ldots, \Si_{n-2}, \ldots, \Si_1, 
\Si_2, \Si_1.$$
Note that if the $\Si_j$ were the generators of $\Cal B_n$, then the element 
$\Delta$ of $\Cal B_n$ determined by the word corresponding to this 
sequence would be a half twist of the entire set of $n$ strings. It is 
easily checked that $\th$ is an $n$-tangle in $B$ whose $j^{th}$ component 
$\th_j$ runs from $a_{0,j}$ to $c_{2m+1,n+1-j}$, and $B_i \cap \th_j$ is 
some $\La_{i,\varphi(i,j)}$. Moreover, given any $j\p > j$, there is some 
$i$ such that $\varphi(i,j\p)=\varphi(i,j)+1$. 

Now let $J_0=\{j_1, \ldots, j_k\}$ be a non-empty subset of 
$\{1, \ldots, n\}$. Let $\thh=\th_{j_1} \cup \cdots \cup \th_{j_k}$. 
Given any \3m $M$ in $B$ such that $\thh \cap M$ is a proper 1-manifold 
in $M$, we denote the exterior of $\thh \cap M$ in $M$ by $M^*$. 

The proof that $\thh$ is excellent is based on a simple strategy which 
may be somewhat obscured by the deluge of notation which follows. 
$B_0$ contains a finite collection of disjoint 3-balls each of which intersects 
$\thh$ in an excellent tangle. As one moves down in $B$ these 3-balls 
expand downwards, following $\thh$ and maintaining the same tangle type. 
As this occurs one notes that the ``empty space''  above and immediately 
around each of these 3-balls consists of disjoint 3-balls each of which meets 
exactly one of these 3-balls in a single disk. They can therefore be 
adjoined to these 3-balls without changing the tangle type. Eventually 
one may encounter a crossing $\Sigma_t$ which involves strings emanating 
from two different  3-balls. At this point the 3-balls come together and 
are attached to a 3-ball just below $\Sigma_t$ whose intersection with 
$\thh$ is an excellent tangle. The result is a new 3-ball whose intersection 
with $\thh$ is, by application of Lemma 6.2, an excellent tangle. One 
then notices that the ``empty space'' in $B$ above and immediately 
around this new 3-ball again can be adjoined to the new 3-ball 
without changing the tangle type. This process 
of downward expansion, amalgamation, and adjunction is then 
continued until it has engulfed all of $B$. 

For $0 \leq i \leq m$ let $J_i$ be the set of those $j \in \{1,\ldots, n\}$ 
such that $\thh \cap B_{i,j} \neq \emptyset$, and let $\Bh_i$ be 
the union of those $B_{i,j}$. Let $T_i=J_0 \cup \cdots \cup J_i$. 
For $1 \leq i \leq m$ let $I_i$ be the set of those $j \in \{1, \ldots, n\}$ 
such that $\thh \cap C_{i,j} \neq \emptyset$, and let $\Ch_i$ be the union 
of those $C_{i,j}$. Let $S_i=I_1 \cup \cdots \cup I_i$. 

We consider how these sets change when $i$ increases by one. We have that 
$\be_{i+1}=\Si_t$ for some $1 \leq t \leq n-1$.

{\it Case 1.} $t$, $t+1 \in J_i$. Then $J_{i+1}=I_{i+1}=J_i$, 
and so $T_{i+1}=T_i$.

{\it Case 2.} $t$, $t+1 \notin J_i$. Again $J_{i+1}=I_{i+1}=J_i$, 
and $T_{i+1}=T_i$.

{\it Case 3.} $t \in J_i$, $t+1 \notin J_i$. Then $I_{i+1}=
J_i \cup \{t+1\}$, $J_{i+1}=I_{i+1}-\{t\}$, and $T_{i+1}=T_i \cup 
\{t+1\}$.

{\it Case 4.} $t \notin J_i$, $t+1 \in J_i$. Then $I_{i+1}=
J_i \cup \{t\}$, $J_{i+1}=I_{i+1}-\{t+1\}$, and $T_{i+1}=T_i \cup 
\{t\}$.

Note that by induction on $i$ it follows that $S_i=T_i$ for 
$1 \leq i \leq m$. It is also easily checked that $T_m=\{1, \ldots, n\}$. 

For $0 \leq i \leq m$ let $R_i$ be the union of all the $B_{r,j}$ and 
$C_{s,j}$ such that $0 \leq r \leq i$, $1 \leq s \leq i$, and $j \in T_i$. 
The components of $R_i$ are, in a sense, the ``minimal'' rectangular 
solids whose union contains $\Bh_0 \cup \Ch_1 \cup \Bh_1 \cup \cdots 
\cup \Ch_i \cup \Bh_i$. Note that $R_0=\Bh_0$ and $R_m=B$. To prove 
the theorem it suffices to prove by induction on $i$ that each component 
of $R^*_i$ is excellent. 

Note that each component $W$ of $\Bh_i$ is a 3-ball which meets $\thh$ in 
$w$ consecutive components of $\La_i$ for some $w \geq 2$. Moreover, there 
is a homeomorphism from $W$ to $B_i$ which is fixed on $\thh \cap W$. 
It follows that $W^*$ is excellent. In particular each component of 
$R^*_0$ is excellent. 

Now suppose the components of $R^*_i$ are excellent. 
Let $P=R_i \cup \Ch_{i+1} \cup \Bh_{i+1}$. 
Let $U$ be the union of all the $B_{r,j}$ and $C_{s,j}$ such that 
$0 \leq r \leq i$, $1 \leq s \leq i$, and $j \in T_{i+1}-T_i$. So $U$ is 
the union of all the blocks of $R_{i+1}$ which lie above a block of $P$ 
but are not contained in $P$. Let $Q=P \cup U$. Let $L$ be the union of all 
the $B_{i+1,j}$ and $C_{i+1,j}$ such that $j \in T_{i+1}-J_{i+1}$. We 
have that $L$ is the union of all the blocks of $R_{i+1}$ which lie below 
a block of $P$ but do not lie in $P$. Then $R_{i+1}=Q \cup L$. We will 
prove in succession that the components of $P^*$, $Q^*$, and 
$R^*_{i+1}$ are excellent. Recall that $\be_{i+1}=\Si_t$. 

{\it Case 1.} $t$, $t+1 \in J_i$. Then $(\Ch_{i+1} \cup \Bh_{i+1})^*$ is 
\homeo\ to $\Bh^*_{i+1}$ and so has all components excellent. The components 
of its intersection with $R^*_i$ are horizontal surfaces with negative 
Euler characteristics whose complements in the boundaries of both 
3-manifolds have no components with closure a disk. Thus by Lemma 6.2 
we have that each component of $P^*$ is excellent. Since $T_{i+1}=T_i$ we 
have $U=\emptyset$ and $L=\emptyset$, and so $R_{i+1}=Q=P$, and we are done. 

{\it Case 2.} $t$, $t+1 \notin J_i$. This case is similar to Case 1. 

{\it Case 3.} $t \in J_i$, $t+1 \notin J_i$. Note that since 
$T_{i+1}-T_i=\{t+1\}$ we have that $U=B_{0,t+1} \cup C_{1,t+1} \cup 
B_{1,t+1} \cup \cdots \cup C_{i,t+1} \cup B_{i,t+1}$ if $i > 0$, and 
$U=B_{0,t+1}$ if $i=0$, and thus is the vertical stack of blocks above 
$C_{i+1,t+1}$. Let $U\p$ denote the union of the corresponding bricks, while 
$U_t=N_{0,t} \cup K_{1,t} \cup N_{1,t} \cup \cdots \cup K_{i,t} \cup N_{i,t}$ if 
$i>0$ and $U_t=N_{0,t}$ if $i=0$, and $U_{t+1}$ is the corresponding 
union with the index $t+1$ in place of $t$. Note that 
$L$ includes the block $B_{i+1,t}$ directly under 
$C_{i+1,t}$. Let $X$ be the component of $R_i$ containing $B_{i,t}$. 
Let $Y$ be the component of $\Ch_{i+1} \cup \Bh_{i+1}$ containing 
$B_{i+1,t+1}$. Denote the components of $P$, $Q$, and $R_{i+1}$ 
containing $X$ by $P_X$, $Q_X$, and $R_{i+1,X}$, with similar notation 
for the components of $P^*$, $Q^*$, and $R^*_{i+1}$ containing $X^*$. 

{\it Subcase (a).} $t+2 \notin J_i$. Then each component of $R_i$ meets a single 
component of $\Ch_{i+1} \cup \Bh_{i+1}$, and vice versa. For components 
of $R_i$ other than $X$ the situation is as in Case 1. Let 
$Z=C_{i+1,t} \cup C_{i+1,t+1} \cup B_{i+1,t+1}$. Then $Z^*$ is \homeo\ 
to $B^*_{i+1,t+1}$ and so is excellent. 

Suppose $t-1 \notin J_i$. Then $Y=Z$,  $P^*_X=X^* \cup Z^*$,  
$Q^*_X=P^*_X \cup U\p \cup U_{t+1}$, and 
$R^*_{i+1,X}=Q^*_X \cup B\p_{i+1,t} \cup N_{i+1,t-1}$. 
We have that $X^* \cap Z^*$ is the surface 
$B^*_{i,t} \cap C^*_{i+1,t}$, which is easily seen to satisfy the 
requirements of Lemma 6.2, and thus $P^*_X$ is excellent. 
$U\p \cup U_{t+1}$ is a 3-ball which meets $P^*_X$ in a disk; thus 
$Q^*_X$ is \homeo\ to $P^*_X$. We have that $B\p_{i+1,t} \cup N_{i+1,t-1}$ 
is a 3-ball which meets $Q^*_X$ in a disk, and so $R^*_{i+1}$ is 
\homeo\ to $Q^*_X$, and we are done. 

Suppose $t-1 \in J_i$. Let $\Yt$ be the closure in $B$ of $Y-Z$. 
Then $P^*_X=X^* \cup \Yt^* \cup Z^* \cup N_{i+1,t-1}$, 
$Q^*_X=P^*_X \cup U\p \cup U_{t+1}$, and 
$R^*_{i+1,X}=Q^*_X \cup B\p_{i+1,t}$. Now $\Yt^*$ is \homeo\ to 
$W^*$, where $W$ is the component of $\Bh_{i+1}$ containing 
$B_{i+1,t-1}$, and so is excellent. $\Yt^*$ meets $X^*$ along a 
surface satisfying the hypotheses of Lemma 6.2, and so $X^* \cup \Yt^*$ 
is excellent. Lemma 6.2 also implies that $X^* \cup \Yt^* \cup Z^*$ is 
excellent. This manifold meets the 3-ball $N_{i+1,t-1}$ in a disk and so 
is \homeo\ to $P^*_X$. As before $Q^*_X$ is \homeo\ to $P^*_X$. The 
3-ball $B\p_{i+1,t}$ meets $Q^*_X$ in a disk, and so $R^*_{i+1,X}$ is 
\homeo\ to $Q^*_X$. 

{\it Subcase (b).}  $t+2 \in J_i$. Let $V$ be the component of $R_i$ 
containing $B_{i,t+2}$. For components of $R_i$ other than $X$ or $V$ 
the situation is as in Case 1. Let $W$ be the component of $\Bh_{i+1}$ 
containing $B_{i+1,t+1}$. So $W=B_{i+1,t+1} \cup \cdots \cup B_{i+1,t+r}$ 
for some $r>1$. Let $Z=C_{i+1,t} \cup W \cup C_{i+1,t+1} \cup \cdots 
\cup C_{i+1,t+r}$. Then $Z^*$ is \homeo\ to $W^*$ and so is excellent. 

Suppose $t-1 \notin J_i$. Then $Y=Z$, $P^*_X=X^* \cup Z^* \cup V^*$, 
$Q^*_X=P^*_X \cup U\p$, and $R^*_{i+1,X}=Q^*_X \cup B\p_{i+1,t} \cup 
N_{i+1,t-1}$. Lemma 6.2 implies that $P^*_X$ is excellent. $U\p$ is a 
3-ball which meets $P^*_X$ in a disk, and so $Q^*_X$ is \homeo\ to $P^*_X$. 
Since $B\p_{i+1,t} \cup N_{i+1,t-1}$ is a 3-ball which meets $Q^*_X$ in 
a disk we have that $R^*_{i+1,X}$ is \homeo\ to $Q^*_X$. 

Suppose $t-1 \in J_i$. Let $\Yt$ be the closure in $B$ of $Y-Z$. Then 
$P^*_X=X^* \cup \Yt^* \cup Z^* \cup V^* \cup N_{i+1,t-1}$, 
$Q^*_X=P^*_X \cup U\p$, and $R^*_{i+1,X}=Q^*_X \cup B\p_{i+1,t}$. 
Successive applications of Lemma 6.2 show that $X^* \cup \Yt^* \cup 
Z^* \cup V^*$ is excellent. Since this manifold meets the 3-ball $N_{i+1,t-1}$ 
in a disk it is \homeo\ to $P^*_X$. For the same reasons $P^*_X$ is 
\homeo\ to $Q^*_X$, which is \homeo\ to $R^*_{i+1,X}$. 

{\it Case 4.} $t \notin J_i$, $t+1 \in J_i$. This case is similar to Case 3. 
\qed \enddemo


\Refs

\ref \no1
\by J. S. Birman
\book Braids, Links, and Mapping Class Groups, {\rm Ann. Math. Studies No. 82}
\publ Princeton University Press
\publaddr Princeton, NJ
\yr 1976
\endref

\ref \no2
\by M. G. Brin and T. L. Thickstun
\paper Open, irreducible 3-manifolds which are end 1-movable
\jour Topology
\vol 26 \issue 2
\pages 211--233
\yr 1987
\endref

\ref \no3
\by M. G. Brin and T. L. Thickstun
\paper 3-manifolds which are end 1-movable
\jour Mem. Amer. Math. Soc. 
\vol 81 (411)
\yr 1989
\endref


\ref \no4
\by M. G. Brin and T. L. Thickstun
\paper Deforming proper homotopy equivalences to 
homeomorphisms in dimension 3
\paperinfo preprint
\endref

\ref \no5
\by E. M. Brown
\paper Contractible 3-manifolds of finite genus at infinity
\jour Trans. Amer. Math. Soc.
\vol 245
\pages 503--514
\yr 1978
\endref

\ref \no6
\by E. M. Brown, M. S. Brown, and C. D. Feustel
\paper On properly embedding planes in 3-manifolds
\jour Proc. Amer. Math. Soc. 
\vol 55
\pages 461--464
\yr 1978
\endref

\ref \no7
\by E. M. Brown and C. D. Feustel
\paper On properly embedding planes in arbitrary 3-manifolds
\jour Proc. Amer. Math. Soc. 
\vol 94
\pages 173--178
\yr 1985
\endref


\ref \no8
\by M. S. Brown
\paper Constructing isotopies on noncompact 3-manifolds
\jour Trans. Amer. Math. Soc. 
\vol 180
\pages 237--263
\yr 1973
\endref

\ref \no9
\by R. H. Fox and E. Artin
\paper Some wild cells and spheres in three-dimensional space
\jour Ann. of Math. (2)
\vol 49
\yr 1948
\pages 979--990
\endref

\ref \no10
\by J. Hempel
\book 3-Manifolds, {\rm Ann. Math. Studies No. 86}
\publ Princeton University Press
\publaddr Princeton, NJ 
\yr 1976
\endref

\ref \no11
\by W. Jaco
\book Lectures on Three-Manifold Topology, {\rm C.B.M.S. Regional Conference 
Series in Mathematics No. 43}
\publ American Mathematical Society
\publaddr Providence, RI 
\yr 1980
\endref

\ref \no12
\by J. W. Morgan
\paper On Thurston's uniformization theorem for three-dimensional manifolds 
\inbook The Smith Conjecture, {\rm Pure and Applied Mathematics Series No. 112} 
\pages 37--126
\publ Academic Press
\publaddr Orlando, FL
\yr 1984
\endref


\ref \no13
\by R. Myers
\paper Homology cobordisms, link concordances, and hyperbolic 3-manifolds
\jour Trans. Amer. Math. Soc. 
\vol 278
\yr 1983
\pages 271--288
\endref

\ref \no14
\by R. Myers
\paper Excellent 1-manifolds in compact 3-manifolds
\jour Topology Appl. 
\vol 49
\yr 1993
\pages 115--127
\endref

\ref \no15
\by R. Myers
\paper End sums of irreducible open 3-manifolds
\paperinfo in preparation
\endref



\ref \no16
\by D. Rolfsen
\book Knots and Links
\publ Publish or Perish Press
\publaddr Berkeley, CA
\yr 1976
\endref



\ref \no17
\by P. Scott and T. Tucker
\paper Some examples of exotic non-compact 3-manifolds
\jour Quart. J. Math. Oxford (2) 
\vol 40
\yr 1989
\pages 481--499
\endref

\ref \no18
\by T. Tucker
\paper On the Fox-Artin sphere and surfaces in non-compact 3-manifolds
\jour Quart. J. Math. Oxford (2) 
\vol 28
\yr 1989
\pages 243--253
\endref



\ref \no19
\by B. N. Winters
\book Proper Planes in Whitehead Manifolds of Finite Genus at Infinity
\bookinfo Ph.D. Thesis 
\publ Oklahoma State University 
\publaddr Stillwater, OK
\yr 1989
\endref

\ref \no20
\by B. N. Winters
\paper Planes in 3-manifolds of finite genus at infinity
\paperinfo preprint
\endref

\ref \no21
\by D. G. Wright
\paper Contractible open manifolds which are not covering spaces 
\jour Topology
\vol 31 \issue 2
\yr 1992
\pages 281-291
\endref


\endRefs
\enddocument